\documentclass[aap,MSNbibl,citesort,dvips]{arximspdf}
\usepackage{mathbh}

\doi{10.1214/11-AAP785}
\volume{22}
\issue{2}
\pubyear{2012}
\firstpage{792}
\lastpage{826}

\makeatletter

\newtheorem{prop}{Proposition}
\newtheorem{theorem}{Theorem}
\newtheorem{lemma}{Lemma}

\newproclaim{defi}{Definition}
\newproclaim{Remark}{Remark}

\newtheorem{corollary}{Corollary}

\newcommand{\C}{{\mathbb C}}
\newcommand{\N}{{\mathbb N}}
\newcommand{\R}{{\mathbb R}}
\newcommand{\Z}{{\mathbb Z}}
\newcommand{\PP}{{\mathbb P}}
\newcommand{\E}{{\mathbb E}}

\newcommand{\cal}{\mathcal}
\newcommand{\eps}{\varepsilon}

\newcommand{\var}{\operatorname{var}}

\makeatother

\begin{document}
\begin{frontmatter}

\title{A scaling analysis of a cat and mouse~Markov~chain}
\runtitle{A scaling analysis of a cat and mouse Markov chain}

\begin{aug}
\author[A]{\fnms{Nelly} \snm{Litvak}\thanksref{t1}\ead[label=e1]{N.Litvak@math.utwente.nl}} and
\author[B]{\fnms{Philippe} \snm{Robert}\corref{}\ead[label=e2]{Philippe.Robert@inria.fr}\ead[label=u1,url]{http://www-rocq.inria.fr/\textasciitilde robert}}
\runauthor{N. Litvak and P. Robert}
\affiliation{University of Twente and INRIA Paris---Rocquencourt}
\address[A]{Faculty of Electrical Engineering\\
Mathematics and Computer Science\\
Department of Applied Mathematics\\
University of Twente\\
7500 AE Enschede\\
The Netherlands\\
\printead{e1}} 
\address[B]{INRIA Paris, Rocquencourt\\
Domaine de Voluceau\\
78153 Le Chesnay\\
France\\
\printead{e2}\\
\printead{u1}}
\end{aug}

\thankstext{t1}{Supported by The Netherlands Organisation
for Scientific Research (NWO) under Meervoud Grant 632.002.401.}

\received{\smonth{10} \syear{2010}}
\revised{\smonth{4} \syear{2011}}

%
\begin{abstract}
If $(C_n)$ is a Markov chain on a discrete state space ${\cal S}$, a
Markov chain $(C_n,M_n)$ on the product space ${\cal S}\times{\cal
S}$, the cat and mouse Markov chain, is constructed. The first
coordinate of this Markov chain behaves like the original Markov chain
and the second component changes only when both coordinates are
equal. The asymptotic properties of this Markov chain are
investigated. A representation of its invariant measure is, in
particular, obtained. When the state space is infinite it is shown
that this Markov chain is in fact null recurrent if the initial
Markov chain $(C_n)$ is positive recurrent and reversible. In
this context, the scaling properties of the location of the
second component, the mouse, are investigated in various
situations: simple random walks in $\Z$ and $\Z^2$ reflected a simple
random walk in $\N$ and also in a continuous time setting. For several
of these processes, a~time scaling with rapid growth gives an
interesting asymptotic behavior related to limiting results for
occupation times and rare events of Markov processes.
\end{abstract}

%

\end{frontmatter}

\section{Introduction}
The PageRank algorithm of Google, as designed by Brin and Page
\cite{BP} in 1998,
describes the web as an oriented graph ${\cal S}$ whose nodes are the
web pages and the
html links between these web pages, the links of the graph. In this
representation, the
importance of a page is defined as its weight for the stationary
distribution of the
associated random walk on the graph. Several off-line algorithms can be
used to estimate
this equilibrium distribution on such a huge state space, they
basically use numerical
procedures (matrix-vector multiplications). See Berkhin
\cite{Berkhin}, for example.
Several on-line algorithms that update the ranking scores while
exploring the graph have
been recently proposed to avoid some of the shortcomings of off-line
algorithms, in
particular, in terms of computational complexity.

The starting point of this paper is an algorithm designed by Abiteboul
et al. \cite{APC} to
compute the stationary distribution of a finite recurrent Markov chain.
In this setting,
to each node of the graph is associated a number, the ``cash'' of the
node. The algorithm
works as follows: at a given time, the node $x$ with the largest value
$V_x$ of cash is
visited, $V_x$ is set to $0$ and the value of the cash of each of its
$d_x$ neighbors is
incremented by $V_x/d_x$. Another possible strategy to update cash
variables is as
follows: a random walker updates the values of the cash at the nodes of
its random path in
the graph. This policy is referred to as the Markovian variant. Both
strategies have the
advantage of simplifying the data structures necessary to manage the
algorithm. It turns
out that the asymptotic distribution, in terms of the number of steps
of the algorithm, of
the vector of the cash variables gives an accurate estimation of the equilibrium
distribution; see Abiteboul et al.~\cite{APC} for the complete
description of the procedure
to get the invariant distribution. See also Litvak and Robert
\cite{Litvak01}. The
present paper does not address the problem of estimating the accuracy
of these
algorithms, it analyzes the asymptotic properties of a simple Markov
chain which appears
naturally in this context.

\subsection*{Cat and mouse Markov chain}
$\!\!$It has been shown in Litvak and Robert~\cite{Litvak01} that, for the
Markovian variant of the algorithm, the distribution of the vector of
the cash variables can be represented with the conditional
distributions of a Markov chain $(C_n,M_n)$ on the discrete state
space ${\cal S}\times{\cal S} $. The sequence~$(C_n)$, representing
the location of the cat, is a Markov chain with transition matrix
$P=(p(x,y))$ associated to the random walk on the graph~${\cal S}$.
The second coordinate, the location of the mouse, $(M_n)$ has the
following dynamic:
\begin{itemize}[--]
\item[--] If $M_n\not=C_n$, then $M_{n+1}=M_n$,
\item[--] If $M_n=C_n$, then, conditionally on $M_n$, the random variable
$M_{n+1}$ has distribution $(p(M_n,y), y\in{\cal S})$ and is
independent of $C_{n+1}$.
\end{itemize}
This can be summarized as follows: the cat moves according to the
transition matrix $P=(p(x,y))$ and the mouse stays idle unless the cat
is at the same site, in which case the mouse also moves independently
according to $P=(p(x,y))$.

The terminology ``cat and mouse problem'' is also used in a somewhat
different way in game theory, the cat playing the role of
the ``adversary.'' See Coppersmith et al. \cite{Copper} and
references therein.

The asymptotic properties of this interesting Markov chain $(C_n,M_n)$
for a~number of transition matrices $P$ are the subject of this paper. In
particular, the asymptotic behavior of the location mouse $(M_n)$ is
investigated. The distribution of $(M_n)$ plays an important role in
the algorithm designed by Abiteboul et al. \cite{APC}; see Litvak and
Robert \cite{Litvak01} for further details. It should be noted that
$(M_n)$ is not, in general, a Markov chain.

\subsection*{Outline of the paper}
Section \ref{secCM} analyzes the recurrence properties of the Markov
chain $(C_n,M_n)$ when the Markov chain $(C_n)$ is recurrent.
A~representation of the invariant measure of $(C_n,M_n)$ in terms of the
reversed process of $(C_n)$ is given.

Since the mouse moves only when the cat arrives at its location, it may seem
quite likely that the mouse will spend most of the time at nodes which
are unlikely for the cat. It is shown that this is indeed the case
when the state space is finite and if the Markov chain $(C_n)$ is
reversible but not in general.

When the state space is infinite and if the Markov chain $(C_n)$ is
reversible, it turns out that the Markov chain $(C_n,M_n)$ is in fact
null recurrent. A precise description of the asymptotic behavior of
the sequence $(M_n)$ is done via a scaling in time and space for
several classes of simple models. Interestingly, the scalings used
are quite diverse, as it will be seen. They are either related to
asymptotics of rare events of ergodic Markov chains or to limiting
results for occupation times of recurrent random walks:
\begin{enumerate}[(1)]
\item[(1)] \textit{Symmetric simple random walks}. The cases of symmetric
simple random walks in $\Z^d$ with $d=1$ and $2$ are analyzed in
Section \ref{secRW}. Note that for $d\geq3$ the Markov chain
$(C_n)$ is transient so that in this case the location of the mouse
does not change with probability $1$ after some random time:
\begin{enumerate}[--]
\item[--] In the one-dimensional case, $d=1$, if $M_0=C_0=0$, on the
linear time scale $t\to n t$, as $n$ gets large, it is shown
that the location of the mouse is of the order of $\sqrt[4]{n}$.
More precisely, the limit in distribution of the process $(M_{\lfloor
nt\rfloor}/\sqrt[4]{n}, t\geq0)$ is a Brownian motion
$(B_1(t))$ taken at the local time at $0$ of another independent
Brownian motion $(B_2(t))$. See Theorem \ref{SRWS} below.

This result can be (roughly) described as follows. Under this linear
time scale the
location of the cat, a simple symmetrical random walk, is of the order
of $\sqrt{n}$ by
Donsker's theorem. It turns out that it will encounter ${\sim}\sqrt
{n}$ times the
mouse. Since the mouse moves only when it encounters the cat and that
it also follows the
sample path of a simple random walk, after $\sqrt{n}$ steps its order
of magnitude will be
therefore of the order of $\sqrt[4]{n}$.
\item[--] When $d=2$, on the linear time scale $t\to n t$, the location of
the mouse is of the order of $\sqrt{\log n}$. More precisely, the
finite marginals of the rescaled processes $(M_{\lfloor
\exp(nt)\rfloor}/\sqrt{n}, t\geq0)$ converge to the
corresponding finite marginals of a Brownian motion in $\R^2$ on a
time scale which is an independent \textit{discontinuous} stochastic
process with independent and nonhomogeneous increments.
\end{enumerate}
\item[(2)] \textit{Reflected simple random walk}. Section \ref{secMM1}
investigates the reflected simple random
walk on the integers. A jump of size
$+1$ (resp., $-1$) occurs with probability $p$ [resp., $(1-p)$]
and the quantity $\rho=p/(1-p)$ is assumed to be strictly less
than $1$ so that the Markov chain $(C_n)$ is ergodic.

If the location of the mouse is far away from the origin,
that is, $M_0=n$ with $n$ large and the cat is at equilibrium, a
standard result shows that it takes a duration of time of the order
of $\rho^{-n}$ for the cat to hit the mouse. This suggests an
exponential time scale $t\to\rho^{-n}t$ to study the evolution of
the successive locations of the mouse. For this time scale it is
shown that the location of the mouse is still of the order of $n$ as
long as $t<W$ where $W$ is some nonintegrable random variable. At
time $t=W$ on the exponential time scale, the mouse has hit $0$ and
after that time the process $(M_{\lfloor
t\rho^{-n}\rfloor}/n)$ oscillates between $0$ and above $1/2$ on
every nonempty time interval.
\item[(3)] \textit{Continuous time random walks}.
Section \ref{secMMI} introduces the cat and mouse process for
continuous time Markov processes. In particular, a discrete
Ornstein--Uhlenbeck process, the $M/M/\infty$ queue, is
analyzed. This is a~birth and death process whose birth rates are
constant and the death rate at $n\in\N$ is proportional to $n$.
When $M_0=n$, contrary to the case of the reflected random walk,
there does not seem to exist a time scale for which a nontrivial
functional theorem holds for the corresponding rescaled
process. Instead, it is possible to describe the asymptotic behavior
of the location of the mouse after the $p$th visit of the cat. It
has a multiplicative representation of the form $nF_1F_2\cdots F_p$
where $(F_p)$ are i.i.d. random variables on $[0,1]$.
\end{enumerate}
The examples analyzed are quite specific. They are, however, sufficiently
representative of the different situations for the dynamic of the
mouse:
\begin{longlist}[(1)]
\item[(1)] One considers the case when an integer valued Markov chain
$(C_n)$ is ergodic and the initial
location of the mouse is far away from $0$.
The correct time scale to investigate the evolution of the location of
the mouse is
given by the duration of time for the occurrence of a rare event for
the original Markov
chain. When the cat hits the mouse at this level, before returning to
the neighborhood
of $0$, it changes the location of the mouse by an additive (resp.,
multiplicative) step
in the case of the reflected random walk (resp., $M/M/\infty$ queue).
\item[(2)] For null recurrent homogeneous random walks, the distribution of
the duration of times between two visits of the cat to the mouse do
not depend on the location of the mouse but it is
nonintegrable. The main problem is therefore to get a functional
renewal theorem associated to an i.i.d. sequence~$(T_n)$ of\vadjust{\goodbreak}
nonnegative random variables such that $\E(T_1)=+\infty$. More
precisely, if
\[
N(t)=\sum_{i\geq1}\mathbh{1}_{\{T_1+\cdots+T_i\leq t\}},
\]
one has to find $\phi(n)$ such that the sequence of processes
$(N(nt)/\phi(n), t\geq0)$ converges as $n$ goes to infinity. When the
tail distribution of $T_1$ has a polynomial decay, several technical
results are available. See Garsia and Lamperti~\cite{Garsia}, for
example. This assumption is nevertheless not valid for the
two-dimensional case. In any case, it turns out that the best way
(especially for $d=2$) to get such results is to formulate the problem
in terms of occupation times of Markov processes for which several
limit theorems are available. This is the key of the results in
Section~\ref{secRW}.
\end{longlist}
The fact that for all the examples considered
jumps occur on the nearest neighbors does not change this qualitative
behavior. Under more general conditions analogous results should hold.
Additionally, this simple setting has the advantage of providing
explicit expressions for most of the constants involved.

\section{The cat and mouse Markov chain}\label{secCM}
In this section we consider a general transition matrix $P=(p(x,y),
x,y\in{\cal S})$ on a discrete state space ${\cal S}$. Throughout the
paper, it is assumed that $P$ is aperiodic, irreducible without loops,
that is, $p(x,x)=0$ for all $x\in{\cal S}$ and with an invariant measure
$\pi$. Note that it is not assumed that $\pi$ has a finite mass. The
sequence $(C_n)$ will denote a Markov chain with transition matrix
$P=(p(x,y))$. It will represent the sequence of nodes which are
sequentially updated by the random walker.

The transition matrix of the reversed Markov chain $(C^*_n)$ is denoted by
\[
p^*(x,y)=\frac{\pi(y)}{\pi(x)} p(y,x)
\]
and, for $y\in{\cal S}$, one defines
\[
H_y^*=\inf\{n>0\dvtx C^*_n=y\} \quad\mbox{and}\quad H_y=\inf\{n>0\dvtx C_n=y\}.
\]
The Markov chain $(C_n,M_n)$ on ${\cal S}\times{\cal S}$ referred to
as the ``cat and mouse Markov chain'' is introduced. Its transition
matrix $Q=(q(\cdot,\cdot))$ is defined as follows: for $x$, $y$,
$z\in{\cal S}$,
%
%
\begin{equation}\label{CMchain}
\cases{
q[(x,y),(z,y)] = p(x,z), &\quad if $x\not=y$;\cr
q[(y,y),(z,w)] = p(y,z)p(y,w).}
\end{equation}
The process $(C_n)$ [resp., $(M_n)$] will be defined as the position of
the cat (resp., the mouse). Note that the position $(C_n)$ of the cat
is indeed a Markov chain with transition matrix $P=(p(\cdot,\cdot))$.
The position of the mouse $(M_n)$ changes only when the cat is at the
same position. In this case, starting from $x\in{\cal S}$ they both
move independently according to the stochastic vector $(p(x,\cdot))$.

Since the transition matrix of $(C_n)$ is assumed to be irreducible
and aperiodic, it is not difficult to check that the Markov chain
$(C_n,M_n)$ is aperiodic and visits with probability $1$ all the
elements of the diagonal of ${\cal S}\times{\cal S}$. In particular,
there is only one irreducible component. Note that $(C_n,M_n)$ itself
is not necessarily irreducible on ${\cal S}\times{\cal S}$, as the
following example shows: take ${\cal S}=\{0,1,2,3\}$ and the
transition matrix $p(0,1)=p(2,3)=p(3,1)=1$ and $p(1,2)=1/2=p(1,0)$; in
this case the element $(0,3)$ cannot be reached starting from $(1,1)$.
\begin{theorem}[(Recurrence)]\label{thInv}
The Markov chain $(C_n,M_n)$ on ${\cal S}\times{\cal S}$ with
transition matrix $Q$
defined by relation (\ref{CMchain}) is recurrent: the measure $\nu$
defined as
%
%
\begin{equation}\label{Rinv}
\nu(x,y)=\pi(x) \E_x\Biggl(\sum_{n=1}^{H_y^*}p(C_n^*,y)\Biggr),\qquad
x, y\in{\cal S},
\end{equation}
is invariant. Its marginal on the second coordinate is given by, for
$y\in{\cal S}$,
\[
\nu_2(y)\stackrel{\mathrm{def.}}{=}\sum_{x\in{\cal S }} \nu(x,y
)=\E_\pi(p(C_0,y)H_y),
\]
and it is equal to $\pi$ on the diagonal, $\nu(x,x)=\pi(x)$ for
$x\in{\cal S}$.
\end{theorem}

In particular, with probability $1$, the elements of ${\cal S}\times
{\cal S}$ for which
$\nu$ is nonzero are visited infinitely often and $\nu$ is, up to a
multiplicative
coefficient, the unique invariant measure. The recurrence property is
not surprising:
the positive recurrence property of the Markov chain $(C_n)$ shows that
cat and mouse
meet infinitely often with probability one. The common location at
these instants is a
Markov chain with transition matrix $P$ and therefore recurrent.
Note that the total mass of $\nu$,
\[
\nu({\cal S}\times{\cal S})=\sum_{y\in{\cal S}} \E_\pi(p(C_0,y)H_y)
\]
can be infinite when ${\cal S}$ is countable.
See Kemeny et al. \cite{Kemeny02} for an introduction on recurrence
properties of
discrete countable Markov chains.

The measure $\nu_2$ on ${\cal S}$ is related to the location of the
mouse under the invariant measure $\nu$.
\begin{pf*}{Proof of Theorem \ref{thInv}}
$\!$From the ergodicity of $(C_n)$ it is clear that~$\nu(x,y)$ is finite
for $x$, $y\in{\cal S}$.
One has first to check that $\nu$ satisfies the equations of invariant
measure for the
Markov chain $(C_n,M_n)$,
%
%
\begin{equation}\label{eqInv}
\nu(x,y)=\sum_{z\not=y}\nu(z,y)p(z,x)\\
+\sum_{z} \nu(z,z)p(z,x)p(z,y),\quad x,y\in{\cal S}.
\end{equation}
For $x$, $y\in{\cal S}$,
%
%
\begin{eqnarray}\label{aaux1}
\sum_{z\not=y}\nu(z,y)p(z,x)
&=&\sum_{z\not=y}\pi(x)p^*(x,z)
\E_z\Biggl(\sum_{n=1}^{H_y^*}p(C_n^*,y)\Biggr)\nonumber\\[-9pt]\\[-9pt]
&=&\pi(x)\E_x\Biggl(\sum_{n=2}^{H_y^*}p(C_n^*,y)\Biggr)
\nonumber
\end{eqnarray}
and
%
%
\begin{eqnarray}
\label{eqaaux3}
&&
\sum_{z\in{\cal S}} \nu(z,z)p(z,x)p(z,y)\nonumber\\[-9pt]\\[-9pt]
&&\qquad=\sum_{z\in{\cal S}}\pi
(x)p^*(x,z)p(z,y)\E_z\Biggl(\sum_{n=0}^{H_z^*-1}p(C_n^*,z)\Biggr).\nonumber
\end{eqnarray}
The classical renewal argument for the invariant distribution $\pi$ of the
Markov chain $(C_n^*)$, and any bounded function $f$ on ${\cal S}$,
gives that
\[
\E_\pi(f)=\frac{1}{\E_z(H_z^*)} \E_z\Biggl(\sum_{n=0}^{H_z^*-1}
f(C_n^*)\Biggr);
\]
see Theorem 3.2, page 12, of Asmussen \cite{Asmussen}, for example.
In particular, we have $\pi(z)=1/\E_z(H_z^*)$, and
%
%
\begin{eqnarray}
\label{eqaaux4}
\E_z\Biggl(\sum_{n=0}^{H_z^*-1}p(C_n^*,z)\Biggr)&=&\E_z(H_z^*) \E
_\pi(p(C_0^*,z))=\frac{\sum_{x\in S}\pi(x)p(x,z)}{\pi(z)}\nonumber\\[-9pt]\\[-9pt]
&=&\frac
{\pi(z)}{\pi(z)}=1.\nonumber
\end{eqnarray}
Substituting the last identity into (\ref{eqaaux3}), we obtain
%
%
\begin{eqnarray}\label{aaux2}
\sum_{z\in{\cal S}} \nu(z,z)p(z,x)p(z,y)&=&\sum_{z\in{\cal S}}\pi
(x)p^*(x,z)p(z,y)\nonumber\\[-9pt]\\[-9pt]
&=&\pi(x)\E_x(p(C_1^*,y)).\nonumber
\end{eqnarray}
Relations (\ref{eqInv})--(\ref{eqaaux3}) and (\ref{aaux2}) show
that $\nu$ is indeed an invariant distribution. At the same time, from
(\ref{eqaaux4}) one gets
the identity $\nu(x,x)=\pi(x)$ for $x\in{\cal S}$.

The second marginal is given by, for $y\in{\cal S}$,
\begin{eqnarray*}
\sum_{x\in{\cal S }} \nu(x,y)
&=&\sum_{t\geq1}\sum_{x\in{\cal S }} \pi(x)\E_x
\bigl(p(C_t^*,y)\mathbh{1}_{\{H_y^*\geq t\}}\bigr)\\[-2pt]
&=&\sum_{t\geq1} \E_\pi\bigl(p(C_t^*,y)\mathbh{1}_{\{
H_y^*\geq t\}}\bigr)\\[-2pt]
&=&\sum_{x\in S}\sum_{z_1,\ldots,z_{t-1}\ne y}\sum_{z_t\in S}\pi
(x)p^*(x,z_1)p^*(z_1,z_2)\cdots p^*(z_{t-1},z_t)p(z_t,y)\\[-2pt]
&=&\sum_{x\in S}\sum_{z_1,\ldots,z_{t-1}\ne y}\sum_{z_t\in
S}p(z_1,x)p(z_2,z_1)\cdots p(z_t,z_{t-1})\pi(z_t)p(z_t,y)\\[-2pt]
&=& \sum_{t\geq1} \E_\pi\bigl(p(C_0,y)\mathbh{1}_{\{H_y\geq
t\}}\bigr)=\E
_\pi(p(C_0,y)H_y),
\end{eqnarray*}
and the theorem is proved.\vspace*{-1pt}
\end{pf*}

The representation (\ref{Rinv}) of the invariant measure can be
obtained (formally) through an
iteration of the equilibrium equations (\ref{eqInv}).
Since the first coordinate of $(C_n,M_n)$ is a Markov chain with
transition matrix $P$
and $\nu$ is the invariant measure for
$(C_n,M_n)$, the first marginal of $\nu$ is thus equal to~$\alpha\pi
$ for some $\alpha>0$, that is,
\[
\sum_{y}\nu(x,y)=\alpha\pi(x),\qquad x\in{\cal S}.
\]
The constant $\alpha$ is in fact the total mass of $\nu$.
In particular, from (\ref{Rinv}), one gets that the quantity
\[
h(x)\stackrel{\mathrm{def.}}{=}\sum_{y\in{\cal S}}\E_x\Biggl(\sum
_{n=1}^{H_y^*}p(C_n^*,y)\Biggr),\qquad x \in{\cal S},
\]
is independent of $x\in{\cal S}$ and equal to $\alpha$. Note that the
parameter
$\alpha$ can be infinite.\vspace*{-1pt}
\begin{prop}[(Location of the mouse in the reversible case)]
If $(C_n)$ is a reversible Markov chain, with the definitions of the
above theorem,
for $y\in{\cal S}$, the relation
\[
\nu_2(y)=1-\pi(y)
\]
holds. If the state space ${\cal S}$ is countable, the Markov chain
$(C_n,M_n)$ is then null recurrent.\vspace*{-1pt}
\end{prop}
\begin{pf}
For $y\in{\cal S}$, by reversibility,
\begin{eqnarray*}
\nu_2(y)&=&\E_\pi(p(C_0,y) H_y)=
\sum_{x} \pi(x)p(x,y)\E_x(H_y)\\[-2pt]
&=&\sum_{x} \pi(y)p(y,x)\E_x(H_y)=\pi(y)\E_y(H_y-1)\\[-2pt]
&=&1-\pi(y).
\end{eqnarray*}
The proposition is proved.\vspace*{-1pt}
\end{pf}
\begin{corollary}[(Finite state space)]\label{corol1}
If the state space ${\cal S}$ is finite with cardinality $N$,\vadjust{\goodbreak} then
$(C_n,M_n)$ converges
in distribution to $(C_\infty,M_\infty)$ such that
%
%
\begin{equation}\label{FinInv}
\PP(C_\infty=x,M_\infty=y)=\alpha^{-1} \pi(x) \E_x\Biggl(\sum
_{n=1}^{H_y^*}p(C_n^*,y)\Biggr) ,\quad x, y\in{\cal S},\vspace*{-1pt}
\end{equation}
with
\[
\alpha=\sum_{y\in{\cal S}}\E_\pi(p(C_0,y)H_y)\vspace*{-1pt}
\]
in particular, $\PP(C_\infty=M_\infty=x)=\alpha^{-1}\pi(x)$. If the
Markov chain $(C_n)$ is reversible, then
\[
\PP(M_\infty=y)=\frac{1-\pi(y)}{N-1}.\vspace*{-1pt}
\]
\end{corollary}

Tetali \cite{Tetali} showed, via linear algebra, that if $(C_n)$ is a
general recurrent Markov chain, then
%
%
\begin{equation}\label{eqconj}
\sum_{y\in{\cal S}} \E_\pi(p(C_0,y) H_y)\le N-1.\vspace*{-1pt}
\end{equation}
See also Aldous and Fill \cite{Aldous}.
It follows that the value $\alpha=N-1$ obtained for reversible chains
is the maximal
possible value of $\alpha$. The constant $\alpha^{-1}$ is the
probability that the cat and mouse are at the same
location.

In the reversible case, Corollary \ref{corol1} implies the intuitive
fact that the less likely a site is for the
cat, the more likely it is for the mouse. This is, however, false in general.
Consider a Markov chain whose state space ${\cal S}$ consists of $r$
cycles with
respective sizes $m_1, \ldots, m_r$ with one common node $0$,
\[
{\cal S}=\{0\} \cup\bigcup_{k=1}^r\{(k,i)\dvtx 1\leq i\leq m_k\},\vspace*{-1pt}
\]
and with the following transitions: for $1\leq k\leq r$ and $2\leq
i\leq m_k$,
\[
p\bigl((k,i),(k,i-1)\bigr)=1,\qquad
p((k,1),0)=1 \quad\mbox{and}\quad p(0,(k,m_k))=\frac{1}{r}.\vspace*{-1pt}
\]
Define $m=m_1+m_2+\cdots+m_r$. It is easy to see that
\[
\pi(0)=\frac{r}{m+r} \quad\mbox{and}\quad \pi(y)=\frac{1}{m+r},\qquad
y\in{\cal S}-\{0\}.\vspace*{-1pt}
\]
One gets that for the location of the mouse, for $y\in{\cal S}$,
\[
\nu_2(y)=E_\pi(p(C_0,y)H_y)=
\cases{
\pi(y)(m-m_k+1), &\quad if $y=(k,m_k), 1\leq k\leq r$,\cr
\pi(y), &\quad otherwise.}\vspace*{-1pt}
\]
Observe that for any $y$ distinct from $0$ and $(k,m_k)$, we have $\pi
(0)>\pi(y)$
and $\nu_2(0)>\nu_2(y)$; the probability to find a mouse in $0$ is
larger than in
$y$. Note that in this example one easily obtains $c=1/r$.\vspace*{-1pt}

\section{Random walks in $\Z$ and $\Z^2$}\label{secRW}
In this section the asymptotic behavior of the mouse when\vadjust{\goodbreak} the cat
follows a recurrent
random walk in $\Z$ and $\Z^2$ is analyzed. The jumps of the cat are uniformly
distributed on the neighbors of the current location.

\subsection{One-dimensional random walk}\label{ssec1}
The transition matrix $P$ of this random walk is given by
\[
p(x,x+1)= \tfrac{1}{2}=p(x,x-1),\qquad x\in\Z.
\]
\subsection*{Decomposition into cycles}
If the cat and the mouse start at the same location, they stay together a
random duration of time $G$ which is geometrically distributed with
parameter $1/2$. Once they are at different locations for the first
time, they are at distance $2$ so that the duration of time $T_2$ until
they meet again has the same distribution as the hitting time of $0$
by the random walk which starts at $2$. The process
\[
\bigl((C_n,M_n), 0\leq n\leq G+T_2\bigr)
\]
is defined as a \textit{cycle}. The sample path of the Markov chain $(C_n,M_n)$
can thus be decomposed into a sequence of cycles. It should be noted
that, during a~cycle, the mouse moves only during the period with
duration $G$.

Since one investigates the asymptotic properties of the sample paths
of~$(M_n)$ on the linear time scale $t\to nt$ for $n$ large,
to get limit theorems one should thus estimate the number of cycles
that occur in a time interval $[0,nt]$. For this purpose, we compare
the cycles of the cat and mouse process to the cycles of a simple
symmetric random walk, which are the time intervals between two
successive visits to zero by the process $(C_n)$.
Observe that a~cycle of $(C_n)$ is equal to $1+T_1$, where $T_1$ is the
time needed to reach zero starting from 1. Further, $T_2$ is the sum of
two independent random variables distributed as $T_1$. Hence, one
guesses that on the linear time scale $t\to nt$ the number of cycles on
$[0,nt]$ for $(C_n,M_n)$ is asymptotically equivalent to $1/2$ of the
number of cycles on $[0,nt]$ for $(C_n)$, as $n\to\infty$. It is well
known that the latter number is of the order $\sqrt{n}$. Then the
mouse makes order of~$\sqrt{n}$ steps of a simple symmetric random
walk, and thus its location must be of the order $\sqrt[4]{n}$.

To make this argument precise, we first prove technical Lemma \ref
{lem1}, which says that only $o(\sqrt{n})$ of $(C_n)$-cycles can be
fitted into the time interval of the order $\sqrt{n}$. Next,
Lemma \ref{lem2} proves that the number of cycles of length $T_2+2$ on
$[0,nt]$, scaled by $\sqrt{n}$, converges to $1/2$ of the local time of
a Brownian motion, analogously to the corresponding result for the
number of cycles of a simple symmetric random walk \cite{Knight}.
Finally, the main limiting result for the location of the mouse is
given by Theorem \ref{SRWS}.
\begin{lemma}\label{lem1}
For any $x$, $\eps>0$ and $K>0$,
\[
\lim_{n\to+\infty}\PP\Biggl(\inf_{0\leq k\leq\lfloor x\sqrt
{n}\rfloor}\frac{1}{\sqrt{n}}
\sum_{i=k}^{k+\lfloor\eps\sqrt{n}\rfloor}(1+T_{1,i})\leq K\Biggr)=0,
\]
where $(T_{1,i})$ are i.i.d. random variables with the same
distribution as the first
hitting time of $0$ of $(C_n)$, $T_1=\inf\{n>0\dvtx C_n=0$
with $C_0=1\}$.
\end{lemma}
\begin{pf}
If $E$ is\vspace*{1pt} an exponential random variable with parameter $1$ independent of
the sequence $(T_{1,i})$, by using the fact that, for $u\in(0,1)$,
$\E(u^{T_1})=(1-\sqrt{1-u^2})/u$, then for $n\geq2$,
\[
\log\PP\Biggl(\frac{1}{\sqrt{n}}\sum_{i=0}^{\lfloor\eps\sqrt
{n}\rfloor} (1+T_{1,i})\leq E\Biggr)=
\lfloor\eps\sqrt{n}\rfloor\log\bigl(1-\sqrt{1-e^{-2/\sqrt
{n}}}\bigr)
\leq-\eps\sqrt[4]{n}.
\]
Denote by
\[
m_n=\inf_{0\leq k\leq\lfloor x\sqrt{n}\rfloor}\frac{1}{\sqrt{n}}
\sum_{i=k}^{k+\lfloor\eps\sqrt{n}\rfloor}(1+T_{1,i})
\]
the above relation gives
\[
\PP(m_n\leq E)\leq\sum_{k=0}^{\lfloor x\sqrt{n}\rfloor
} \PP\Biggl(\frac{1}{\sqrt{n}}\sum_{i=k}^{k+\lfloor\eps\sqrt
{n}\rfloor} (1+T_{1,i})\leq E\Biggr)\leq\bigl(\bigl\lfloor x\sqrt{n}\bigr\rfloor
+1\bigr) e^{-\eps\sqrt[4]{n}},
\]
hence,
\[
\sum_{n=2}^{+\infty} \PP(m_n\leq E)<+\infty
\]
and, consequently, with probability $1$, there exists $N_0$ such that,
for any $n\geq N_0$,
we have $m_n> E$. Since $\PP(E\geq K)>0$, the lemma is proved.
\end{pf}
\begin{lemma}\label{lem2}
Let, for $n\geq1$, $(T_{2,i})$ i.i.d. random variables with the same
distribution
as $T_2=\inf\{k>0\dvtx C_k=0$ with $C_0=2\}$ and
\[
u_n=\sum_{\ell=1}^{+\infty} \mathbh{1}_{\{\sum_{k=1}^{\ell}
(2+T_{2,k})< n\}},
\]
then the process
$(u_{\lfloor t n \rfloor}/\sqrt{n} )$
converges in distribution to $(L_B(t)/2)$,\break where~$L_B(t)$ is the local
time process at time $t\geq0$
of a standard Brownian motion.
\end{lemma}
\begin{pf}
The variable $T_2$ can be written as a sum $T_1+T_1'$ of independent
random variables
$T_1$ and $T_1'$ having the same distribution as $T_1$ defined in the
above lemma. For
$k\geq1$, the variable $T_{2,k}$ can be written as
$T_{1,2k-1}+T_{1,2k}$. Clearly,
\[
\frac{1}{2} \sum_{\ell=1}^{+\infty} \mathbh{1}_{\{\sum
_{k=1}^{\ell}(1+T_{1,k})< n\}} -\frac{1}{2} \leq u_n
\leq
\frac{1}{2} \sum_{\ell=1}^{+\infty} \mathbh{1}_{\{\sum
_{k=1}^{\ell} (1+T_{1,k})< n\}}.
\]
Furthermore,
\[
\Biggl(\sum_{\ell=1}^{+\infty} \mathbh{1}_{\{\sum
_{k=1}^{\ell} (1+T_{1,k})< n\}}, n\geq1\Biggr)\stackrel
{\mathrm{dist.}}{=}
(r_n)\stackrel{\mathrm{def.}}{=}\Biggl(\sum_{\ell=1}^{n-1} \mathbh
{1}_{\{C_\ell=0\}}, n\geq1\Biggr),
\]
where $(C_n)$ is the symmetric simple random walk.

A classical result by Knight \cite{Knight} (see also Borodin
\cite{Borodin} and
Perkins \cite{Perkins}) gives that the process $(r_{\lfloor n t\rfloor
}/\sqrt{n})$ converges in distribution to
$(L_B(t))$ as $n$ gets large. The lemma is proved.
\end{pf}

The main result of this section can now be stated.
\begin{theorem}[(Scaling of the location of the mouse)] \label{SRWS}
If $(C_0,M_0)\in\N^2$, the convergence in distribution
\[
\lim_{n\to+\infty}\biggl(\frac{1}{\sqrt[4]{n}}M_{\lfloor nt\rfloor
}, t\geq0\biggr)
\stackrel{\mathrm{dist.}}{=}\bigl(B_1(L_{B_2}(t)),
t\geq0\bigr)
\]
holds, where $(B_1(t))$ and $(B_2(t))$ are independent standard
Brownian motions on $\R$ and
$(L_{B_2}(t))$ is the local time process of $(B_2(t))$ at $0$.
\end{theorem}

The location of the mouse at time $T$ is therefore of the
order of $\sqrt[4]{T}$ as~$T$ gets large. The limiting process can be
expressed as a Brownian motion slowed down by the process of the local
time at $0$ of an independent Brownian motion. The quantity
$L_{B_2}(T)$ can be interpreted as the scaled duration of time the cat
and the mouse spend together.
\begin{pf*}{Proof of Theorem \ref{SRWS}}
Without loss of generality, one can assume that \mbox{$C_0=M_0$}.
A coupling argument is used. Take:
\begin{itemize}[--]
\item[--] i.i.d. geometric random variables $(G_i)$ such that $\PP
(G_1\geq p)=1/2^{p-1}$ for \mbox{$p\geq1$};
\item[--] $(C_{k}^a)$ and $(C_{j,k}^b)$, $j\geq1$, i.i.d. independent
symmetric random walks
starting from $0$;
\end{itemize}
and assume that all these random variables are independent. One
denotes, for $m=1$,
$2$ and $j\geq1$,
\[
T^{b}_{m,j}=\inf\{k\geq0\dvtx C_{j,k}^b=m\}.
\]
Define
\[
(C_k,M_k)=
\cases{(C_k^a,C_k^a), &\quad $0\leq k< G_1$,\cr
(C_{G_1}^a-2I_1+I_1C_{1,k-G_1}^b,C_{G_1}^a), &\quad $G_1\leq k\leq\tau_1$,}
\]
with $I_1=C^a_{G_1}-C^a_{G_1-1}$, $\tau_1=G_1+T^{b}_{2,1}$.
It is not difficult to check that
\[
[(C_k,M_k), 0\leq k\leq\tau_1]
\]
has the same
distribution as the cat and mouse Markov chain during a~cycle as
defined above.

Define $t_0=0$ and $t_i=t_{i-1}+\tau_i$, $s_0=0$ and $s_{i}=s_{i-1}+G_{i}$.
The $(i+1)$th cycle is defined as
\[
(C_k,M_k)=
\cases{(C_{k-t_i+s_i}^a,C_{k-t_i+s_i}^a), \qquad t_i\leq k<
t_i+G_{i+1},\cr
(C_{s_{i+1}}^a-2I_{i+1}+I_{i+1}C_{i+1,k-t_i-G_{i+1}}^b,C_{s_{i+1}}^a),\vspace*{2pt}\cr
\hspace*{103.59pt}\quad
t_i+G_{i+1}\leq k\leq t_{i+1},}
\]
with $I_{i+1}=C^a_{s_{i+1}}-C^a_{s_{i+1}-1}$ and $\tau
_{i+1}=G_{i+1}+T^{b}_{2,i+1}$.
The sequence $(C_n,M_n)$ has the same distribution as the Markov chain with
transition matrix $Q$ defined by relation (\ref{CMchain}).

With this representation, the location $M_n$ of the mouse at time $n$
is given by~$C^a_{\kappa_n}$, where $\kappa_n$ is the number of steps the mouse
has made up to time $n$, formally defined as
\[
\kappa_n\stackrel{\mathrm{def.}}{=}\sum_{i=1}^{+\infty}
\Biggl[\sum_{\ell=1}^{i-1} G_\ell+(n-t_{i-1})\Biggr]\mathbh{1}_{\{
t_{i-1}\leq n\leq t_{i-1}+G_{i}\}}
+\sum_{i=1}^{+\infty}
\Biggl[\sum_{\ell=1}^{i} G_\ell\Biggr] \mathbh{1}_{\{
t_{i-1}+G_{i}<n<t_{i}\}},
\]
in particular,
%
%
\begin{equation}\label{aux1}
\sum_{\ell=1}^{\nu_n} G_\ell\leq\kappa_n\leq\sum_{\ell=1}^{\nu
_n+1} G_\ell
\end{equation}
with $\nu_n$ defined as the number of cycles of the cat and mouse
process up to time $n$:
\[
\nu_n=\inf\{\ell\dvtx t_{\ell+1}> n\}=\inf\Biggl\{\ell\dvtx \sum
_{k=1}^{\ell+1}(G_k+T_{2,k}^{b})> n\Biggr\}.
\]
Define
\[
\overline{\nu}_n=\inf\Biggl\{\ell\dvtx \sum_{k=1}^{\ell+1}
(2+T_{2,k}^{b})> n\Biggr\},
\]
then, for $\delta>0$, on the event
$\{\overline{\nu}_n> \nu_n+\delta\sqrt{n}\}$,
\begin{eqnarray*}
n&\geq&\sum_{k=1}^{\nu_n+\delta\sqrt{n}}(2+ T_{2,k}^{b}
) \geq\sum_{k=1}^{\nu_n+1}
[G_k+T_{2,k}^{b}]+\sum_{k=\nu_n+2}^{\nu_n+\delta\sqrt
{n}}(2+T_{2,k}^{b})
-\sum_{k=1}^{\nu_n+1}(G_k-2)\\
&\geq& n+\sum_{k=\nu_n+2}^{\nu_n+\delta\sqrt{n}}
(2+T_{2,k}^{b})
-\sum_{k=1}^{\nu_n+1}(G_k-2).
\end{eqnarray*}
Hence,
\[
\sum_{k=\nu_n+2}^{\nu_n+\delta\sqrt{n}}(2+T_{2,k}^{b}
)\leq
\sum_{k=1}^{\nu_n+1}(G_k-2);\vadjust{\goodbreak}
\]
since $T_{1,k}^b\leq2+T_{2,k}^b$, the relation
%
%
\begin{eqnarray}\label{aux2}
\bigl\{\overline{\nu}_n- \nu_n> \delta
\sqrt{n}\bigr\}&\subset&\Biggl\{\inf_{1\leq
\ell\leq\nu_n}\sum_{k=\ell}^{\ell+\lfloor\delta\sqrt{n}\rfloor
} T_{1,k}^b\leq
\sum_{k=1}^{\nu_n+1}(G_k-2), \overline{\nu}_n>\nu_n\Biggr\}\hspace*{-30pt}
\nonumber\\[-8pt]\\[-8pt]
&\subset&\Biggl\{\inf_{1\leq
\ell\leq\overline{\nu}_n}\sum_{k=\ell}^{\ell+\lfloor\delta
\sqrt{n}\rfloor} T_{1,k}^b\leq
\sup_{1\leq\ell\leq\overline{\nu}_n}\sum_{k=1}^{\ell
+1}(G_k-2)\Biggr\}\nonumber
\end{eqnarray}
holds. Since $\E(G_1)=2$, Donsker's theorem gives the following
convergence in distribution:
\[
\lim_{K\to+\infty}
\Biggl( \frac{1}{\sqrt{K}}\sum_{k=1}^{\lfloor t K\rfloor+1}(G_k-2),
0\leq t\leq1\Biggr)
\stackrel{\mathrm{dist.}}{=} \bigl(\var(G_1)W(t), 0\leq t\leq
1\bigr),
\]
where $(W(t))$ is a standard Brownian motion, and, therefore,
%
%
\begin{equation}\label{egg}
\lim_{K\to+\infty}\frac{1}{\sqrt{K}}\sup_{1\leq\ell\leq K}\sum
_{k=1}^{\ell+1}(G_k-2)
\stackrel{\mathrm{dist.}}{=} \var(G_1)\sup_{0\leq t\leq1} W(t).
\end{equation}
For $t>0$, define
\[
\bigl(\Delta_n(s), 0\leq s\leq t\bigr)\stackrel{\mathrm{def.}}{=}
\biggl(\frac{1}{\sqrt{n}} \bigl(\overline{\nu}_{\lfloor n s\rfloor}-\nu
_{\lfloor n s\rfloor}\bigr),0\leq s\leq t\biggr).
\]
By relation (\ref{aux2}) one gets that, for $0\leq s\leq t$,
%
%
\begin{eqnarray}\label{ehh}
\{ \Delta_n(s) >\delta\}&\subset&\Biggl\{\inf_{1\leq
\ell\leq\overline{\nu}_{\lfloor ns\rfloor}}\sum_{k=\ell}^{\ell
+\lfloor
\delta\sqrt{n}\rfloor} T_{1,k}^b\leq
\sup_{1\leq\ell\leq\overline{\nu}_{\lfloor ns\rfloor}}\sum
_{k=1}^{\ell+1}(G_k-2)\Biggr\}\nonumber\\[-8pt]\\[-8pt]
&\subset&\Biggl\{\inf_{1\leq\ell\leq\overline{\nu}_{\lfloor
nt\rfloor}}\sum_{k=\ell}^{\ell+\lfloor\delta\sqrt{n}\rfloor}
T_{1,k}^b\leq\sup_{1\leq\ell\leq\overline{\nu}_{\lfloor
nt\rfloor}}\sum_{k=1}^{\ell+1}(G_k-2)\Biggr\}.
\nonumber
\end{eqnarray}
Letting $\eps>0$, by Lemma \ref{lem2} and relation (\ref{egg}), there
exist some $x_0>0$ and~$n_0$ such that if $n\geq n_0$, then, respectively,
%
%
\begin{equation}\label{eii}\qquad
\PP\bigl(\overline{\nu}_{\lfloor n t\rfloor}\geq x_0\sqrt
{n}
\bigr)\leq\eps\quad\mbox{and}\quad
\PP\Biggl(\sup_{1\leq\ell\leq x_0\sqrt{n} }\sum_{k=1}^{\ell
+1}(G_k-2) \geq x_0\sqrt{n}\Biggr)\leq\eps.
\end{equation}
By using relation (\ref{ehh}),
\begin{eqnarray*}
&&\Bigl\{\sup_{0\leq s\leq t}\Delta_n(s) > \delta\Bigr\}\\
&&\qquad\subset\bigl\{\overline{\nu}_{\lfloor n t\rfloor}\geq x_0\sqrt{n}\bigr\}\\
&&\qquad\quad{}\cup
\Biggl\{\inf_{1\leq\ell\leq\overline{\nu}_{\lfloor nt\rfloor
}}\sum_{k=\ell}^{\ell+\lfloor\delta\sqrt{n}\rfloor}
T_{1,k}^b\leq
\sup_{1\leq\ell\leq\overline{\nu}_{\lfloor
nt\rfloor}}\sum_{k=1}^{\ell+1}(G_k-2), \overline{\nu}_{\lfloor n
t\rfloor}<x_0\sqrt{n} \Biggr\}\\
&&\qquad\subset\bigl\{\overline{\nu}_{\lfloor n t\rfloor}\geq x_0\sqrt{n}\bigr\}
\cup
\Biggl\{\inf_{1\leq\ell\leq x_0\sqrt{n}}\sum_{k=\ell}^{\ell
+\lfloor\delta\sqrt{n}\rfloor} T_{1,k}^b\leq
\sup_{1\leq\ell\leq x_0\sqrt{n}}\sum_{k=1}^{\ell+1}(G_k-2)
\Biggr\}.
\end{eqnarray*}
With a similar decomposition with the partial sums of $(G_k-2)$,
relations~(\ref{eii}) give the inequality, for $n\geq n_0$,
\[
\PP\Bigl(\sup_{0\leq s\leq t}\Delta_n(s) > \delta\Bigr)
\leq2\eps+
\PP\Biggl(\inf_{1\leq k\leq x_0\sqrt{n}}\frac{1}{\sqrt{n}}
\sum_{i=k}^{k+\lfloor\delta\sqrt{n}\rfloor} T_{1,i}^b\leq
x_0\Biggr).
\]
By Lemma \ref{lem1}, the left-hand side is thus arbitrarily small if
$n$ is
sufficiently large. In a similar way the same results holds for the variable
$\sup(-\Delta_n(s)\dvtx\allowbreak {0\leq s\leq t})$.
The variable $\sup(|\Delta_n(s)|\dvtx 0\leq s\leq t)$ converges therefore in
distribution to $0$. Consequently, by using
relation (\ref{aux1}) and the law of large numbers, the same property
holds for
\[
\sup_{0\leq s\leq t}\frac{1}{\sqrt{n}}\bigl(\kappa_{\lfloor n
s\rfloor}-2\overline{\nu}_{\lfloor n s\rfloor}\bigr).
\]
Donsker's theorem\vspace*{1pt} gives that the sequence of processes \mbox{$(C^a_{\lfloor
\sqrt{n} s\rfloor}/\sqrt[4]{n},
0\leq s\leq t)$} converges in distribution to $(B_1(s), 0\leq s\leq t)$.
In particular, for
$\eps$ and $\delta>0$, there exists some $n_0$ such that if $n\geq
n_0$, then
\[
\PP\biggl(\sup_{0\leq u,v\leq t, |u-v|\leq\delta} \frac{1}{\sqrt[4]{n}}
\bigl|C^a_{\lfloor\sqrt{n} u\rfloor}-C^a_{\lfloor\sqrt{n}
v\rfloor}\bigr|\geq\delta\biggr)\leq\eps;
\]
see Billingsley \cite{Billingsley02}, for example.
Since $M_n=C_{\kappa_n}^a$ for any $n\geq1$, the processes
\[
\biggl(\frac{1}{\sqrt[4]{n}}M_{\lfloor n s\rfloor}, 0\leq s\leq
t\biggr) \quad\mbox{and}\quad
\biggl(\frac{1}{\sqrt[4]{n}}C_{2\overline{\nu}_{\lfloor n s\rfloor
}}^a, 0\leq s\leq t\biggr)
\]
have therefore the same asymptotic behavior for the convergence in distribution.
Since, by construction $(C_k^a)$ and $(\overline{\nu}_n)$ are
independent, with Skorohod's
representation theorem, one can assume that, on an appropriate
probability space with two
independent Brownian motions $(B_1(s))$ and $(B_2(s))$, the convergences
\begin{eqnarray*}
\lim_{n\to+\infty}\bigl(C^a_{\lfloor\sqrt{n} s\rfloor}/\sqrt[4]{n},
0\leq s\leq t\bigr)&=&\bigl(B_1(s),
0\leq s\leq t\bigr),\\
\lim_{n\to+\infty} \bigl(\overline{\nu}_{\lfloor n s\rfloor}/\sqrt
{n}\bigr)&=&\bigl(L_{B_2}(s)/2, 0\leq s\leq t\bigr)
\end{eqnarray*}
hold almost surely for the norm of the supremum. This concludes the
proof of the theorem.
\end{pf*}

\subsection{Random walk in the plane}\label{ssec2}
The transition matrix $P$ of this random walk is given by, for $x\in\Z^2$,
\[
p\bigl(x,x+(1,0)\bigr) =p\bigl(x,x-(1,0)\bigr)=p\bigl(x,x+(0,1)\bigr)
=p\bigl(x,x-(0,1)\bigr)= \tfrac{1}{4}.
\]

\subsection*{Decomposition into cycles}
In the one-dimensional case, when the cat and the mouse start at the
same location, when they are separated for the first time, they are at
distance $2$, so that the next meeting time has the same distribution
as the hitting time of $0$ for the simple random walk when it starts
at $2$. For $d=2$, because of the geometry, the situation is more
complicated. When the cat and the mouse are separated for the first
time, there are several possibilities for the patterns of their
respective locations and not only one as for $d=1$. A finite Markov
chain has to be introduced that describes the relative position of the
mouse with respect to the cat.

Let $e_1=(1,0)$, $e_{-1}=-e_1$, $e_2=(0,1)$, $e_{-2}=-e_2$
and the set of
unit vectors of $\Z^2$ is denoted by ${\cal E}=\{e_1,e_{-1},
e_2,e_{-2}\}$.
Clearly, when the cat and the mouse are at the same site, they stay
together a geometric number of
steps whose mean is~$4/3$. When they are just separated, up to a
translation, a symmetry or a rotation, if the mouse is at $e_1$, the
cat will be at $e_2$,
$e_{-2}$ or $-e_1$ with probability $1/3$. The next time the cat will
meet the mouse
corresponds to one of the instants of visit to ${\cal E}$ by the sequence
$(C_n)$. If one considers only these visits, then, up to a translation,
it is not difficult
to see that the position of the cat and of the mouse is a Markov chain
with transition
matrix $Q_R$ defined below.
\begin{defi}\label{def1}
Let $e_1=(1,0)$, $e_{-1}=-e_1$, $e_2=(0,1)$, $e_{-2}=-e_2$
and the set of
unit vectors of $\Z^2$ is denoted by ${\cal E}=\{e_1,e_{-1},
e_2,e_{-2}\}$.

If $(C_n)$ is a random walk in the plane, $(R_n)$ denotes the sequence
in ${\cal E}$ such
that $(R_n)$ is the sequence of unit vectors visited by $(C_n)$ and
%
%
\begin{equation}\label{Dir}
r_{ef}\stackrel{\mathrm{def.}}{=}\PP(R_1=f\mid R_0=e),\qquad
e, f\in {\cal E}.
\end{equation}
A transition matrix $Q_R$ on ${\cal E}^2$ is defined as follows: for
$e$, $f$, $g\in{\cal E}$,
%
%
\begin{equation}\label{CMR}
\cases{Q_R((e,g), (f,g)) = r_{ef},\qquad e\not=g,\cr
Q_R((e,e), (e,-e)) = 1/3,\cr
Q_R((e,e), (e,\overline{e}))=Q_R\bigl((e,e), (e,-\overline{e})\bigr)=1/3,}
\end{equation}
with the convention that $\overline{e}$, $-\overline{e}$ are the unit
vectors orthogonal
to $e$, $\mu_R$ denotes the invariant probability distribution
associated to $Q_R$ and
$D_{\cal E}$ is the diagonal of ${\cal E}^2$.
\end{defi}

A characterization of the matrix $R$ is as follows. Let
\[
\tau^+=\inf(n>0\dvtx C_n\in{\cal E}) \quad\mbox{and}\quad \tau=\inf(n\geq0\dvtx
C_n\in{\cal E}),
\]
then clearly $r_{ef}=\PP(C_{\tau^+}=f\mid C_0=e)$. For $x\in\Z^2$, define
\[
\phi(x)=\PP(C_{\tau}=e_1\mid C_0=x).
\]
By symmetry, it is easily seen that the coefficients of $R$ can be
determined by $\phi$.
For $x\notin{\cal E}$, by looking at the state of the Markov chain
at time
$1$, one gets the relation
\[
\Delta\phi(x)\stackrel{\mathrm{def.}}{=}\phi(x+e_1)+\phi
(x+e_{-1})+\phi(x+e_2)+\phi(x+e_{-2})-4\phi(x)= 0
\]
and $\phi(e_i)=0$ if $i\in\{-1,2,-2\}$ and $\phi(e_1)=1$. In other
words, $\phi$ is the
solution of a \textit{discrete Dirichlet problem}: it is a harmonic
function (for the discrete
Laplacian) on $\Z^2$ with fixed values on ${\cal E}$. Classically,
there is
a unique solution to the Dirichlet problem; see Norris \cite{Norris},
for example. An
explicit expression of $\phi$ is, apparently, not available.
\begin{theorem}\label{theo2}
If $(C_0,M_0)\in\N^2$, the convergence in distribution of finite marginals
\[
\lim_{n\to+\infty} \biggl(\frac{1}{\sqrt{n}} M_{\lfloor
e^{nt}\rfloor}, t\geq
0\biggr)\stackrel{\mathrm{dist.}}{=} (W(Z(t)))
\]
holds, with
\[
(Z(t))=\biggl(\frac{16\mu_R(D_{\cal E})}{3\pi}L_B(T_t)\biggr),
\]
where $\mu_R$ is the probability distribution on ${\cal E}^2$
introduced in
Definition \ref{def1}, the process $(W(t))=(W_1(t),W_2(t))$ is a
two-dimensional Brownian motion and:
\begin{itemize}[--]
\item[--] $(L_B(t))$ the local time at $0$ of a standard Brownian motion
$(B(t))$ on $\R$
independent of $(W(t))$.
\item[--] For $t\geq0$, $T_t=\inf\{s\geq0\dvtx B(s)=t\}$.
\end{itemize}
\end{theorem}
\begin{pf}
The proof follows the same lines as before: a convenient construction
of the process to
decouple the time scale of the visits of the cat and the motion of the
mouse. The
arguments which are similar to the ones used in the proof of the
one-dimensional case are
not repeated.

Let $(R_n,S_n)$ be the Markov chain with transition matrix $Q_R$ that
describes the relative positions of the cat and the mouse at the
instances of visits of $(C_n,M_n)$ to ${\cal E}\times{\cal E}$ up to
rotation, symmetry and translation. For~$N$ visits to the set
${\cal E}\times{\cal E}$, the proportion of time the cat and the mouse
will have
met is given by
\[
\frac{1}{N}\sum_{\ell=1}^N \mathbh{1}_{\{R_\ell=S_\ell
\}};
\]
this quantity converges almost surely to $\mu_R(D_{{\cal E}})$.\vadjust{\goodbreak}

Now one has to estimate the number of visits of the cat to the set
${\cal E}$.
Kasahara \cite{Kas1} (see also Bingham \cite{Bing2} and
Kasahara \cite{Kas2}) gives that,
for the convergence in distribution of the finite marginals, the
following convergence holds:
\[
\lim_{n\to+\infty} \Biggl(\frac{1}{n} \sum_{i=0}^{\lfloor
e^{nt}\rfloor}\mathbh{1}_{\{C_i\in{\cal E}\}}
\Biggr)\stackrel{\mathrm
{dist.}}{=} \biggl(\frac{4}{\pi}L_B(T_t)\biggr).
\]
The rest of the proof follows the same lines as in the proof of
Theorem \ref{SRWS}.~%
\end{pf}
\begin{Remark*}
Tanaka's Formula (see Rogers and Williams \cite{Rogers06}) gives the relation
\[
L(T_t)=t-\int_0^{T_t}\operatorname{sgn}(B(s)) \,dB(s),
\]
where sgn$(x)=-1$ if $x<0$ and $+1$ otherwise. Since the process
$(T_t)$ has independent
increments and that the $T_t$'s are stopping times, one gets that
$(L(T_t))$ has also
independent increments. With the function $t\to T_t$ being
discontinuous, the limiting process
$(W(Z(t)))$ is also discontinuous. This is related to the fact that the
convergence of
processes in the theorem is minimal: it is only for the convergence in
distribution of finite
marginals. For $t\geq0$, the distribution of $L(T_t)$ is an
exponential distribution with
mean~$2t$; see Borodin and Salminen \cite{BorSal}, for example. The
characteristic function of
\[
W_1\biggl(\frac{16\mu_R(D_{{\cal E}})L(T_t)}{3\pi}\biggr)
\]
at $\xi\in\C$ such that $\operatorname{Re}(\xi)=0$ can be easily
obtained as
\[
\E\bigl(e^{i\xi W_1[Z(t)]}\bigr)=\frac{\alpha_0^2}{\alpha
_0^2+\xi^2t}\qquad\mbox{with }
\alpha_0=\frac{\sqrt{3\pi}}{4\sqrt{\mu_R(D_{{\cal E}})}}.
\]
With a simple inversion, one gets that the density of this random
variable is a bilateral
exponential distribution given by
\[
\frac{\alpha_0}{2\sqrt{t}}\exp\biggl(-\frac{\alpha_0}{\sqrt
{t}}|y|\biggr),\qquad y\in\R.
\]
The characteristic function can be also represented as
\[
\E\bigl(e^{i\xi W_1[Z(t)]}\bigr)=\frac{\alpha_0^2}{\alpha
_0^2+\xi t}
=\exp\biggl(\int_{-\infty}^{+\infty} (e^{i\xi u}-1)\Pi
(t,u) \,du\biggr)
\]
with
\[
\Pi(t,u)=\frac{e^{-\alpha_0|u|/\sqrt{t}}}{|u|},\qquad u\in\R.
\]
$\Pi(t,u) \,du$ is in fact the associated L\'evy measure of the
nonhomogeneous process with
independent increments $(W_1(Z(t)))$. See Chapter 5 of Gikhman and
Skorohod \cite{Gikhman}.\vadjust{\goodbreak}
\end{Remark*}

\section{The reflected random walk}\label{secMM1}
In this section the cat follows a simple ergodic random walk on the
integers with a
reflection at $0$; an asymptotic analysis of the evolution of the
sample paths of the
mouse is carried out. Despite being a quite simple example, it
exhibits already an interesting scaling behavior.

Let $P$ denote the transition matrix of the simple reflected random
walk on $\N$,
%
%
\begin{equation}\label{RW}
\cases{p(x,x+1)= p, &\quad $x\geq0$,\cr
p(x,x-1) = 1-p, &\quad $x\not=0$,\cr
p(0,0)=1-p.}
\end{equation}
It is assumed that $p\in(0,1/2)$ so that the corresponding Markov
chain is positive
recurrent and reversible and its invariant probability distribution is
a geometric random
variable with parameter $\rho\stackrel{\mathrm{def.}}{=}p/(1-p)$. In
this case, one can check that the measure $\nu$ on $\N^2$ defined in
Theorem \ref{thInv} is given by
\[
\cases{\nu(x,y) = \rho^{x}(1-\rho), &\quad $0\leq x<y-1$,\vspace*{1pt}\cr
\nu(y-1,y) = \rho^{y-1}(1-\rho)(1-p) , \vspace*{1pt}\cr
\nu(y,y) = \rho^{y}(1-\rho), \vspace*{1pt}\cr
\nu(y+1,y) = \rho^{y+1}(1-\rho)p, \vspace*{1pt}\cr
\nu(x,y) = \rho^{x}(1-\rho), &\quad $x>y+1$.}
\]
The following proposition describes the scaling for the dynamics of the cat.
\begin{prop}\label{ConvExp}
If, for $n\geq1$, $T_n=\inf\{k>0 \dvtx C_k=n\}$, then, as $n$ goes
to infinity, the
random variable $T_n/\E_0(T_n)$ converges in distribution to an
exponentially distributed
random variable with parameter $1$ and
\[
\lim_{n\to+\infty}\E_0(T_n)\rho^n=\frac{1+\rho}{(1-\rho)^2}
\]
with $\rho={p}/{(1-p)}$.

If $C_0=n$, then $T_0/n$ converges almost surely to $(1+\rho)/(1-\rho)$.
\end{prop}
\begin{pf}
The first convergence result is standard; see Keilson
\cite{Keilson01} for
closely related results. Note that the Markov chain $(C_n)$ has the
same distribution as the embedded Markov chain of the $M/M/1$ queue
with arrival rate $p$ and service rate~$q$. The first part of the
proposition is therefore a~discrete analogue of the convergence
result of Proposition 5.11 of Robert~\cite{Robert08}.

If $C_0=n$ and define by induction $\tau_{n}=0$ and, for $0\leq i\leq n$,
\[
\tau_i=\inf\{k\geq0\dvtx C_{k+\tau_{i+1}}=i\},
\]
hence, $\tau_n+\cdots+\tau_i$ is the first time when the cat crosses
level $i$.
The strong Markov property gives that the $(\tau_i,0\leq i\leq n-1)$\vadjust{\goodbreak}
are i.i.d. A standard calculation (see Grimmett and Stirzaker
\cite{Grimmett09}, e.g.) gives that
\[
\E(u^{\tau_1})=\frac{1-\sqrt{1-4pqu}}{2pu},\qquad 0\leq
u\leq1,
\]
hence, $\E(\tau_0)=(1+\rho)/(1-\rho)$. Since $T_0=\tau
_{n-1}+\cdots+\tau_{0}$, the last part of the proposition is
therefore a consequence of the law of large numbers.\vspace*{-3pt}
\end{pf}

\subsection*{Additive jumps}
An intuitive picture of the main phenomenon is as follows. It is
assumed that the mouse is at level $n$ for some $n$ large. If the cat
starts at~$0$, according to the above proposition, it will take of the
order of $\rho^{-n}$ steps to reach the mouse. The cat and the mouse
will then interact for a short amount of time until the cat returns in
the neighborhood of $0$, leaving the mouse at some new location $M$. Note
that, because $n$ is large, the reflection condition does not play a
role for the dynamics of the mouse at this level and by spatial
homogeneity outside $0$, one has that $M=n+M'$ where $M'$ is some
random variable whose distribution is independent of $n$. Hence, when
the cat has
returned to the mouse $k$ times after hitting $0$ and then went back to
$0$ again, the location of the
mouse can be represented as $n+M'_1+\cdots+M'_k$, where $(M'_i)$ are
i.i.d. with the same distribution as~$M'$. Roughly speaking, on the
exponential time scale $t\to\rho^{-n}t$, it will be seen that the
successive locations of the mouse can be represented with the random
walk associated to $M'$ with a negative drift, that is, $\E(M')<0$.

The
section is organized as follows: one investigates the
properties of the random variable $M'$ and the rest of the section is
devoted to the proof of the functional limit theorem. The main ingredient
is also a decomposition of the sample path of $(C_n,M_n)$ into
cycles. A cycle starts and ends with the cat at $0$ and the
mouse is visited at least once by the cat during the cycle.\vspace*{-3pt}

\subsection*{Free process}
Let $(C'_n,M'_n)$ be the cat and mouse Markov chain associated to the
simple random walk on
$\Z$ without reflection (the free process):
\[
p'(x,x+1)=p=1-p'(x,x-1)\qquad \forall x\in\Z.\vspace*{-1pt}
\]

\begin{prop}\label{PropFree}
If $(C'_0,M'_0)=(0,0)$, then the asymptotic location of the mouse for
the free process
$ M'_\infty=\lim_{n\to\infty} M'_n$ is such that, for
$u\in\C$ such that
$|u|=1$,
%
%
\begin{equation}\label{GenD}
\E(u^{M'_\infty})=\frac{\rho(1-\rho)u^2}{-\rho
^2u^2+(1+\rho)u-1},
\end{equation}
in particular,
\[
\E(M'_\infty)=-\frac{1}{\rho}\quad\mbox{and}\quad\E\biggl(\frac{1}{\rho
^{M'_\infty}}\biggr)=1.
\]
Furthermore, the relation
%
%
\begin{equation}\label{InSq}
\E\biggl(\sup_{n\geq0} \frac{1}{\sqrt{\rho}^{M_n'}}
\biggr)<+\infty\vadjust{\goodbreak}
\end{equation}
holds. If $(S_k)$ is the random walk associated to a sequence of i.i.d. random
variables with the same distribution as $M_\infty'$ and $(E_i)$ are
i.i.d. exponential
random variables with parameter $(1+\rho)/(1-\rho)^2$, then the
random variable~$W$
defined by
%
%
\begin{equation}\label{W}
W=\sum_{k=0}^{+\infty} \rho^{-S_k} E_k
\end{equation}
is almost surely finite with infinite expectation.
\end{prop}
\begin{pf}
Let $\tau=\inf\{n\geq1\dvtx C_n'<M_n'\}$, then, by looking at
the different cases, one has
\[
M_\tau'=
\cases{1, &\quad if $M_1=1, C_1=-1$, \cr
1+M_\tau'', &\quad if $M_1=1, C_1=1$, \cr
-1+M_\tau'', &\quad if $M_1=-1$,}
\]
where $M_\tau''$ is an independent r.v. with the same distribution as
$M_\tau'$.
Hence, for $u\in\C$ such that $|u|=1$, one gets that
\[
\E(u^{M'_\tau})=\biggl((1-p)\frac{1}{u}+p^2u\biggr)\E
(u^{M'_\tau})+p(1-p)u
\]
holds.
Since $M'_\tau-C'_\tau=2$, after time $\tau$, the cat and the mouse
meet again with
probability $\rho^2$. Consequently,
\[
M'_\infty\stackrel{\mathrm{dist.}}{=}\sum_{i=1}^{1+G} M'_{\tau,i},
\]
where $(M'_{\tau,i})$ are i.i.d. random variables with the same
distribution as $M'_\tau$
and $G$ is an independent geometrically distributed random variable
with parameter
$\rho^2$. This identity gives directly the expression (\ref{GenD})
for the characteristic function of $M_\infty'$ and also the relation
$\E(M_\infty')=-1/\rho$.

Recall that the mouse can move one step up only when it is at the same
location as the cat, hence, one gets the upper bound
\[
\sup_{n\geq0} M_n'\leq U\stackrel{\mathrm{def.}}{=}1+ \sup_{n\geq
0} C_n'
\]
and the fact that $U-1$ has the same distribution as the invariant
distribution of the
reflected random walk $(C_n)$, that is, a geometric distribution with
parameter $\rho$ gives
directly inequality (\ref{InSq}).

Let $N=(N_t)$ be a Poisson process with rate $(1-\rho)^2/(1+\rho)$,
then one can check
the following identity for the distributions:
%
%
\begin{equation}\label{ExpInt}
W\stackrel{\mathrm{dist.}}{=}\int_0^{+\infty} \rho^{-S_{N_t}}\, dt.
\end{equation}
By the law of large numbers, $(S_{N_t}/t)$ converges almost surely to
$-(1+\rho)/[(1-\rho)^2\rho]$. One gets\vspace*{1pt} therefore that $W$ is almost
surely finite.
From (\ref{GenD}), one gets $u\mapsto\E(u^{M_\infty'})$ can
be analytically extended to the interval
\[
\frac{1+\rho-\sqrt{(1-\rho)(1+3\rho)}}{2\rho^2}<u<\frac{1+\rho
+\sqrt{(1-\rho)(1+3\rho)}}{2\rho^2}
\]
in particular, for $u=1/\rho$ and
its value is $\E(\rho^{-M_\infty'})=1$. This gives by (\ref{W}) and
Fubini's theorem that $\E(W)=+\infty$.
\end{pf}

Note that $\E(\rho^{-M_\infty'})=1$ implies that the exponential moment
$\E(u^{M'_\infty})$ of the random variable $M'_\infty$
is finite for $u$ in the
interval $[1,1/\rho]$.

\subsection*{Exponential functionals}
The representation (\ref{ExpInt}) shows that the variable~$W$ is an
exponential functional
of a compound Poisson process. See Yor \cite{Yor03}. It can be seen
as the invariant
distribution of the auto-regressive process~$(X_n)$ defined as
\[
X_{n+1}\stackrel{\mathrm{def.}}{=}\rho^{-A_n}X_n+E_n,\qquad n\geq0.
\]
The distributions of these random variables are investigated in Guillemin
et al. \cite{Guillemin05} when $(A_n)$ are nonnegative. See also
Bertoin and Yor \cite{BY}. The above proposition shows that $W$ has a
heavy-tailed distribution. As it will be seen
in the scaling result below, this has a qualitative impact on the
asymptotic behavior of
the location of the mouse. See Goldie \cite{Goldie} for an analysis of
the asymptotic
behavior of tail distributions of these random variables.

\subsection*{A scaling for the location of the mouse}
The rest of the section is devoted to the analysis of the location of
the mouse when it is
initially far away from the location of the cat. Define
\[
s_1=\inf\{\ell\geq0\dvtx C_\ell=M_\ell\} \quad\mbox{and}\quad t_1=\inf\{\ell
\geq s_1\dvtx
C_\ell=0\}
\]
and, for $k\geq1$,
%
%
\begin{equation}\label{tk}\qquad
s_{k+1}=\inf\{\ell\geq t_k\dvtx C_\ell=M_\ell\} \quad\mbox{and}\quad
t_{k+1}=\inf\{\ell\geq s_{k+1}\dvtx
C_\ell=0\}.
\end{equation}
Proposition \ref{ConvExp} suggests an exponential time scale for a
convenient scaling of
the location of the mouse. When the mouse is initially at $n$ and the
cat at the origin, it
takes the duration $s_1$ of the order of $\rho^{-n}$ so that the cat
reaches this
level. Just after that time, the two processes behave like
the free process on~$\Z$ analyzed above, hence, when the cat returns
to the origin (at time
$t_1$), the mouse is at position $n+M_\infty'$. Note that on the
extremely fast exponential time scale $t\to\rho^{-n} t$,
the (finite) time that the cat and mouse spend together is vanishing,
and so is the time needed for the cat to reach zero from \mbox{$n+M_\infty'$}
(linear in $n$ by the second statement of Proposition \ref{ConvExp}).
Hence, on the exponential time scale, $s_1$ is a finite exponential
random variable, and $s_2$ is distributed as a sum of two i.i.d. copies
of~$s_1$. The following proposition presents a precise formulation of
this description, in particular, a proof of the corresponding scaling
results. For the sake of simplicity, and because of the topological
intricacies of
convergence in distribution, in a first step the convergence
result is restricted on the time interval $[0,s_2]$, that is, on the
two first
``cycles.'' Theorem \ref{TheoMM1} below gives the full statement of
the scaling result.
\begin{prop}\label{propconv}
If $M_0=n\geq1$ and $C_0=0$, then, as $n$ goes to infinity, the random variable
$(M_{t_1}-n,\rho^n t_1)$ converges in distribution to $(M_\infty',
E_1)$ and the process
\[
\biggl(\frac{M_{\lfloor t\rho^{-n}\rfloor}}{n}\mathbh{1}_{\{
0\leq t< \rho^{n}s_2\}}\biggr)
\]
converges in distribution for the Skorohod topology to the process
\[
\bigl(\mathbh{1}_{\{t<E_1+\rho^{-M_\infty'}E_2\}
}\bigr),
\]
where the distribution of $M_\infty'$ is as defined in
Proposition \ref{PropFree}, and it is independent of $E_1$ and $E_2$,
two independent
exponential random variables with parameter $(1+\rho)/(1-\rho)^2$.
\end{prop}
\begin{pf}
For $T>0$, ${\cal D}([0,T],\R)$ denotes the space of cadlag functions,
that is, of right
continuous functions with left limits, and $d^0$ is the metric on this
space defined by,
for $x$, $y\in{\cal D}([0,T],\R)$,
\[
d^0(x,y)= \inf_{\varphi\in{\cal H}} \biggl[ \sup_{0\leq s<t<T}
\biggl|\log\frac{\varphi(t)-\varphi(s)}{t-s}\biggr|+
\sup_{0\leq s<T} |x(\varphi(s))-y(s)| \biggr],
\]
where ${\cal H}$ is the set of nondecreasing functions $\varphi$ such
that $\varphi(0)=0$
and \mbox{$\varphi(T)=T$}. See Billingsley \cite{Billingsley02}.

An upper index $n$ is added on the variables $s_1$, $s_2$, $t_1$ to
stress the dependence
on $n$. Take three independent Markov chains $(C^a_k)$, $(C^b_k)$ and
$(C^c_k)$ with
transition matrix $P$ such that $C^a_0=C^c_0=0$, $C^b_0=n$ and, for
$i=a$, $b$, $c$,
$T_p^i$ denotes the hitting time of $p\geq0$ for $(C^i_k)$.
Since\vspace*{-1pt}
$((C_k,M_k), s_1^n\leq
k\leq t_1^n)$ has the same distribution as $((n+C_k',n+M_k'), 0\leq k
<T_0^b)$, by the strong
Markov property, the sequence $(M_k, k\leq s_2^n)$ has the same
distribution as $(N_k,
0\leq k\leq T^a_n+T^b_0+T_n^c)$, where
%
%
\begin{equation}
\label{Nk}
N_k=
\cases{n, &\quad $k\leq T^a_n$,\cr
n+M'_{k-T^a_n}, &\quad $T^a_n\leq k\leq T^a_n+T^b_0$,\cr
n+M'_{T^b_0}, &\quad $T^a_n+T^b_0 \leq k\leq
T^a_n+T^b_0+T_{n+M'_{T^b_0}}^c$.}
\end{equation}
Here $((C_k^b-n,M'_k), 0\leq k\leq T_0^b)$ is a sequence with the same
distribution as the
free process with initial starting point $(0,0)$ and killed at the
hitting time of $-n$ by the first coordinate. Additionally, it
is independent of the Markov chains $(C^a_k)$ and $(C^c_k)$.
In particular, the random variable $M_{t_1}-n$, the jump of the
mouse from its initial position when the cat hits $0$, has the same
distribution as $M'_{T_0^b}$. Since $T_0^b$ converges almost surely to
infinity, $M'_{T_0^b}$ is converging in distribution to $M_\infty'$.

Proposition \ref{ConvExp} and the independence of $(C^a_k)$ and
$(C^c_k)$ show that the
sequences $(\rho^n T_n^a)$ and $(\rho^nT_n^c)$ converge in
distribution to two
independent exponential random variables $E_1$ and $E_2$ with parameter
$(1+\rho)/(1-\rho)^2$. By using Skorohod's Representation theorem, (see
Billingsley \cite{Billingsley02}) up to a~change of probability
space, it can be
assumed that these convergences hold for the almost sure convergence.

By representation (\ref{Nk}), the rescaled process $((M_{\lfloor t\rho
^{-n}\rfloor}/n)\mathbh{1}_{\{0\leq t< \rho^{n}s_2\}},\break
t\leq T)$
has the same distribution as
\[
x_n(t)\stackrel{\mathrm{def.}}{=}
\cases{1, &\quad $t< \rho^n T^a_n$,\vspace*{1pt}\cr
\displaystyle  1+\frac{1}{n}M'_{\lfloor\rho^{-n}t-T^a_n\rfloor}, &\quad
$\rho^n
T^a_n\leq t< \rho^n(T^a_n+T^b_0)$,\vspace*{2pt}\cr
\displaystyle  1+\frac{1}{n}M'_{T^b_0}, &\quad $\rho^n
(T^a_n+T^b_0) \leq t< \rho^n
(T^a_n+T^b_0+T_{n+M'_{T^b_0}}^c)$, \vspace*{2pt}\cr
\displaystyle  0, &\quad $t\geq\rho^n (T^a_n+T^b_0+T_{n+M'_{T^b_0}}^c)$,}
\]
for $t\leq T$.
Proposition \ref{ConvExp} shows that $T_0^b/n$ converges almost surely to
$(1-\rho)/\allowbreak(1+\rho)$ so that $(\rho^n(T^a_n+T^b_0))$ converges to
$E_1$ and, for
$n\geq1$,
\[
\rho^nT_{n+M'_{T^b_0}}^c=\rho^{-M'_{T^b_0}}\rho
^{n+M'_{T^b_0}}T_{n+M'_{T^b_0}}^c
\longrightarrow\rho^{-M'_{\infty}} E_2,
\]
almost surely as $n$ goes to infinity. Additionally, one has also
\[
\lim_{n\to+\infty}\frac{1}{n}\sup_{k\geq0} |M_k'|=0,
\]
almost surely. Define
\[
x_\infty=\bigl(\mathbh{1}_{\{t<T\wedge(E_1+\rho^{-M_\infty
'}E_2)\}}\bigr),
\]
where $a\wedge b=\min(a,b)$ for $a$, $b\in\R$.

Time change. For $n\geq1$ and $t>0$, define $u_n$ (resp., $v_n$) as the minimum
(resp., maximum) of $t\wedge
\rho^n\lceil T^a_n+T^b_0+T_{n+M'_{T^b_0}}^c\rceil$ and $t\wedge
(E_1+\rho^{-M'_{\infty}} E_2)$, and
\[
\varphi_n(s)=
\cases{ \displaystyle  \frac{v_n}{u_n} s, &\quad $0\leq s\leq u_n$,\vspace*{2pt}\cr
\displaystyle  v_n+(s-u_n)\frac{T-v_n}{T-u_n}, &\quad $u_n<s\leq T$.}
\]
Noting that $\varphi_n\in{\cal H}$, by using this function in the
definition of the
distance~$d^0$ on ${\cal D}([0,T],\R)$ to have an upper bound of
$(d(x_n,x_\infty))$
and with the above convergence results, one gets that, almost
surely, the sequence $(d(x_n,x_\infty))$ converges to~$0$. The
proposition is proved.
\end{pf}
\begin{theorem}[(Scaling for the location of the mouse)]\label{TheoMM1}
If $M_0=n$, $C_0=0$, then the process
\[
\biggl(\frac{M_{\lfloor t\rho^{-n}\rfloor}}{n}\mathbh{1}_{\{
t<\rho^n t_n\}}\biggr)
\]
converges in distribution for the Skorohod topology to the process
$(\mathbh{1}_{\{t<W\}})$, where $W$ is the random
variable defined by
(\ref{W}).

If $H_0$ is the hitting time of $0$ by $(M_n)$,
\[
H_0=\inf\{s\geq0\dvtx M_s=0\},
\]
then, as $n$ goes to infinity, $\rho^n H_0$ converges in distribution
to $W$.
\end{theorem}
\begin{pf}
In the same way as in the proof of Proposition \ref{propconv}, it can
be proved that for
$p\geq1$, the random vector $[(M_{t_k}-n,\rho^n t_k), 1\leq k\leq p]$
converges in distribution to
the vector
\[
\Biggl(S_k,\sum_{i=0}^{k-1} \rho^{-S_i}E_i\Biggr)
\]
and, for $k\geq0$, the convergence in distribution
%
%
\begin{equation}\label{Conv1}
\lim_{n\to+\infty} \biggl(\frac{M_{\lfloor t\rho^{-n}\rfloor
}}{n}\mathbh{1}_{\{0\leq t< \rho^{n}t_k\}}\biggr)
=
\bigl(\mathbh{1}_{\{t<E_1+\rho^{-S_1}E_2+\cdots+\rho
^{-S_{k-1}}E_k\}}\bigr)
\end{equation}
holds for the Skorohod topology.

Let $\phi\dvtx [0,1]\to\R_+$ be defined by $\phi(s)=\E(\rho^{-s
M_\infty'})$, then $\phi(0)=\phi(1)=1$
and $\phi'(0)<0$, since $\phi$ is strictly convex then for all $s<1$,
$\phi(s)<1$.

If $C_0=M_0=n$, and the sample path of $(M_k-n, k\geq0)$ follows the sample
path of a reflected random walk starting at $0$, we have, in
particular, that the
supremum of its successive values is integrable. By
Proposition \ref{PropFree}, as $n$ goes to infinity, $M_{t_1}-n$ is
converging in distribution to $M'_{\infty}$. Lebesgue's theorem gives
therefore that the averages are also converging, hence, since~$\E
(M_\infty')$ is negative, there exists
$N_0$ such that if $n\geq N_0$,
%
%
\begin{equation}\label{ineq1}\quad
\E_{(n,n)}(M_{t_1})\stackrel{\mathrm{def.}}{=}\E(M_{t_1}\mid
M_0=C_0=n)\leq
n+\frac{1}{2}\E(M'_{\infty})=n-\frac{1}{2\rho}.
\end{equation}
Note that $t_1$ has the same distribution as $T_0$ in
Proposition \ref{ConvExp} when $C_0=n$. Proposition \ref{ConvExp} now
implies that there exists $K_0\geq0$ so that,
for $n\geq N_0$,
%
%
\begin{equation}\label{ineq2}
\rho^{n/2}\E_{(0,n)}\bigl(\sqrt{t_1} \bigr)\leq K_0.
\end{equation}
The identity $\E(1/\rho^{M_\infty'})=1$ implies that
$\E(\rho^{-M_\infty'/2})<1$, and inequality~(\ref{InSq}) and
Lebesgue's theorem imply that one can choose $0<\delta<1$ and $N_0$,
so that
%
%
\begin{equation}\label{ineq3}
\E\bigl(\rho^{(n-M_{t_1})/2}\bigr)\leq\delta
\end{equation}
holds for $n\geq N_0$.
Let $\nu=\inf\{k\geq1\dvtx M_{t_k}\leq N_0\}$ and, for $k\geq1$, ${\cal
G}_k$ the
$\sigma$-field generated by the random variables $(C_j,M_j)$ for
$j\leq t_k$.
Because of inequality~(\ref{ineq1}), one can check that the sequence
\[
\biggl(M_{t_{k\wedge\nu}}+ \frac{1}{2\rho}(k\wedge\nu), k\geq
0\biggr)
\]
is a super-martingale with respect to the filtration $({\cal G}_k)$, hence,
\[
\E(M_{t_{k\wedge\nu}})+ \frac{1}{2\rho}\E(k\wedge
\nu)\leq\E(M_0)=n.
\]
Since the location of the mouse is nonnegative, by letting $k$ go to
infinity, one gets
that $\E(\nu)\leq2\rho n$. In particular, $\nu$ is almost surely a
finite random
variable.

Intuitively, $t_{\nu}$ is the time when the mouse reaches the area
below a finite boundary $N_0$.
Our goal now is to prove that the sequence $(\rho^n t_{\nu})$
converges in distribution to $W$.
For $p\geq1$ and on the event $\{\nu\geq p\}$,
%
%
\begin{equation}\label{ineq5}
\bigl(\rho^n (t_{\nu}-t_p)\bigr)^{1/2}=\Biggl(\sum_{k=p}^{\nu-1}
\rho^n(t_{k+1}-t_k)\Biggr)^{1/2}
\leq\sum_{k=p}^{\nu-1} \sqrt{\rho^n(t_{k+1}-t_k)}.
\end{equation}
For $k\geq p$, inequality (\ref{ineq2}) and the strong Markov property
give that the relation
\[
\rho^{M_{t_k}/2} \E\bigl[\sqrt{t_{k+1}-t_k}\mid{\cal G}_k\bigr]=
\rho^{M_{t_k}/2} \E_{(0,M_{t_k})}\bigl[\sqrt{t_1}\bigr]\leq K_0
\]
holds on the event $\{\nu>k\}\subset\{M_{t_k}>N_0\}$. One gets
therefore that
\begin{eqnarray*}
\E\bigl(\sqrt{\rho^n(t_{k+1}-t_k)}\mathbh{1}_{\{k<\nu
\}}\bigr)&=&
\E\bigl(\rho^{(n-M_{t_k})/2}\mathbh{1}_{\{k<\nu\}}
\rho^{M_{t_k}/2}
\E\bigl[\sqrt{t_{k+1}-t_k}\mid{\cal G}_k\bigr]\bigr)\\
&\leq& K_0\E\bigl(\rho^{(n-M_{t_k})/2}\mathbh{1}_{\{k<\nu
\}}\bigr)
\end{eqnarray*}
holds, and, with inequality (\ref{ineq3}) and again the strong Markov property,
\begin{eqnarray*}
\E\bigl(\rho^{(n-M_{t_k})/2}\mathbh{1}_{\{k<\nu\}
}\bigr)&=&
\E\bigl(\rho^{-\sum_{j=0}^{k-1}(M_{t_{j+1}}-M_{t_j})/2}\mathbh
{1}_{\{k<\nu\}}\bigr)\\
&\leq&\delta\E\bigl(\rho^{-\sum
_{j=0}^{k-2}(M_{t_{j+1}}-M_{t_j})/2}\mathbh{1}_{\{k-1<\nu
\}}\bigr)
\leq\delta^k.
\end{eqnarray*}
Relation (\ref{ineq5}) gives therefore that
\[
\E\bigl(\sqrt{\rho^n(t_{\nu}-t_p)}\bigr)\leq\frac{K_0\delta
^p}{1-\delta}.
\]
For $\xi\geq0$,
%
%
\begin{eqnarray}\label{Joblot}\quad
|\E(e^{-\xi\rho^nt_{\nu}})-\E(e^{-\xi
\rho^nt_{p}})|
&\leq&\bigl|\E\bigl(1-e^{-\xi\rho^n(t_{\nu}-t_p)^+}\bigr)
\bigr|+\PP(\nu<p)\nonumber\\
&=&\int_0^{+\infty} \xi e^{-\xi u}\PP\bigl( \rho^n(t_{\nu}-t_p)\geq u\bigr)\,
du+\PP(\nu<p)\\
&\leq&\frac{K_0\delta^p}{1-\delta}\int_0^{+\infty} \frac{\xi
}{\sqrt{u}} e^{-\xi u} \,du
+\PP(\nu<p)\nonumber
\end{eqnarray}
by using Markov's inequality for the random variable $\sqrt{\rho
^n(t_\nu-t_p)}$.\break Since~$\rho^n t_p$ converges in distribution to
$E_0+\rho^{-S_1}E_1+\cdots+\rho^{-S_p}E_p$, one can prove that, for
$\eps>0$, by choosing a fixed $p$ sufficiently large and that if $n$
is large enough, then
the Laplace transforms at $\xi\geq0$ of the random variables $\rho
^nt_{\nu}$ and $W$ are
at a distance less than $\eps$.

At time $t_\nu$ the location $M_{t_\nu}$ of the mouse is $x\leq N_0$
and the cat is at~$0$. Since the sites visited by $M_n$ are a Markov chain with
transition matrix
$(p(x,y))$, with probability $1$, the number $R$ of jumps for the mouse
to reach~$0$ is finite. By recurrence of $(C_n)$, almost surely, the cat will
meet the mouse $R$ times in
a finite time. Consequently, if $H_0$ is the time when the mouse hits
$0$ for the first time, then by the strong Markov property, the
difference $H_0-t_\nu$ is
almost surely a finite random variable. The convergence in distribution
of~$(\rho^nH_0)$
to $W$ is therefore proved.
\end{pf}

\subsection*{Nonconvergence of scaled process after $W$}

Theorem \ref{TheoMM1} could suggest that the convergence holds for a
whole time axis, that is,
\[
\lim_{n\to+\infty}\biggl(\frac{M_{\lfloor t\rho^{-n}\rfloor}}{n},
t\geq0\biggr)=\bigl(\mathbh{1}_{\{t<W\}},t\geq0\bigr)
\]
for the Skorohod topology. That is, after time $W$ the rescaled process
stays at $0$ like
for fluid limits of stable stochastic systems. However, it turns out
that this
convergence does not hold at all for the following intuitive (and
nonrigorous) reason. Each time the cat
meets the mouse at $x$ large, the location of the mouse is at
$x+M_\infty'$ when
the cat returns to $0$, where $M_\infty'$ is the random variable
defined in
Proposition \ref{PropFree}. In this way, after the $k$th visit of the
cat, the mouse is at
the $k$th position of a random walk associated to~$M_\infty'$ starting
at $x$. Since
$\E(1/\rho^{M_\infty'})=1$, Kingman's result (see Kingman~\cite{Kingman01}) implies that
the hitting time of $\delta n$, with $0<\delta<1$, by this random walk
started at $0$ is of the order of $\rho^{-\delta n}$.
For each of the steps of the random walk, the cat needs also of the
order of
$\rho^{-\delta n}$ units of time. Hence, the mouse reaches the level
$\delta n$ in order of $\rho^{-2\delta n}$
steps, and this happens on any finite interval $[s,t]$ on the time
scale $t\to\rho^{-n}t$ only if
$\delta\leq1/2$. Thus, it is very likely that the next relation holds:
\[
\lim_{n\to+\infty} \PP\biggl(\sup_{s\leq u\leq t} \frac
{M_{\lfloor
u\rho^{-n}\rfloor}}{n}=\frac{1}{2}\biggr)=1.
\]
Note that this implies that for $\delta\leq1/2$ on the time scale
$t\to\rho^{-n}t$ the mouse will cross the level $\delta n$ infinitely
often on any finite interval! The difficulty in proving this statement
is that the mouse is not at $x+M_\infty'$ when the cat
returns at $0$ at time $\tau_x$ but at $x+M_{\tau_x}'$, so that the
associated random
walk is not space-homogeneous but only asymptotically close to the one
described above. Since an
exponentially large number of steps of the random walks are considered,
controlling the
accuracy of the approximation turns out to be a~problem. Nevertheless,
a partial result is established in the next proposition.
\begin{prop}\label{osc}
If $M_0=C_0=0$, then for any $s$, $t>0$ with $s<t$, the relation
%
%
\begin{equation}\label{Irancy}
\lim_{n\to+\infty} \PP\biggl(\sup_{s\leq u\leq t} \frac
{M_{\lfloor
u\rho^{-n}\rfloor}}{n}\geq\frac{1}{2}\biggr)=1
\end{equation}
holds.
\end{prop}

It should be kept in mind that, since $(C_n,M_n)$ is
recurrent, the process~$(M_n)$ returns infinitely often to $0$ so that
relation (\ref{Irancy}) implies that the scaled process exhibits
oscillations for the norm of the supremum on compact intervals.
\begin{pf*}{Proof of Proposition \ref{osc}}
First it is assumed that $s=0$.
If $C_0=0$ and $T_0=\inf\{k>0\dvtx C_k=0\}$, then, in particular,
$\E(T_0)=1/(1-\rho)$. The set ${\cal C}=\{C_0,\ldots,C_{T_0-1}\}$
is a cycle of the
Markov chain, and denote by $B$ its maximal value. The Markov chain can
be decomposed into
independent cycles $({\cal C}_n, n\geq1)$ with the corresponding
values $(T_0^n)$ and
$(B_n)$ for $T_0$ and~$B$. Kingman's result (see Theorem 3.7 of
Robert \cite{Robert08},
e.g.) shows that there exists some constant $K_0$ such that $\PP
(B\geq n)\sim
K_0\rho^{n}$. Taking $0<\delta<1/2$, for $\alpha>0$,
\[
U_n\stackrel{\mathrm{def.}}{=}\rho^{(1-\delta)n} \sum
_{k=1}^{\lfloor\alpha\rho^{-n}\rfloor}
\bigl[\mathbh{1}_{\{B_k\geq\delta n\}}-\PP(B\geq
\delta n)\bigr],
\]
then, by Chebyshev's inequality, for $\eps>0$,
\begin{eqnarray*}
\PP(|U_n|\geq\eps)&\leq&\rho^{(2-2\delta)n}{\alpha\rho^{-n}}\frac
{\operatorname{Var}(\mathbh{1}_{\{B\geq\delta n\}})}{\eps^2}\le
\frac{\alpha}{\eps^2}\rho^{(1-2\delta)n}\PP(B\geq\delta n)\\
&\leq&
\frac{\alpha K_0}{\eps^2} \rho^{(1-\delta)n}.
\end{eqnarray*}
By using Borel--Cantelli's lemma, one gets that the sequence $(U_n)$
converges almost surely to $0$, hence, almost surely,
%
%
\begin{equation}\label{eq1}
\lim_{n\to+\infty} \rho^{(1-\delta)n} \sum_{k=1}^{\lfloor\alpha
\rho^{-n}\rfloor}
\mathbh{1}_{\{B_k\geq\delta n\}}=\alpha K_0.
\end{equation}
For $x\in\N$, let $\nu_x$ be the number of cycles up to time $x$,
and the strong law of large numbers gives that, almost
surely,
%
%
\begin{equation}
\label{eqslln}
\lim_{x\to+\infty} \frac{\nu_x}{x}=\lim_{x\to+\infty} \frac
{1}{x}\sum_{k=1}^x\mathbh{1}_{\{C_k=0\}}=1-\rho.
\end{equation}
Denote by $x_n\stackrel{\mathrm{def.}}{=}\lfloor\rho^{-n}t\rfloor$.
For $\alpha_0>0$, the probability that the location of the mouse is
never above level
$\delta n$ on the time interval $(0,x_n]$ is
%
%
\begin{eqnarray}\label{eq2}
&&\PP\Bigl(\sup_{1\leq k\leq\lfloor\rho^{-n}t\rfloor}M_k\leq
\delta
n\Bigr)\nonumber\\[-3pt]
&&\qquad\leq
\PP\Biggl(\sup_{1\leq k\leq\lfloor\rho^{-n}t\rfloor}M_k\leq
\delta
n, \rho^{(1-\delta)n}\sum_{i=0}^{\nu_{x_n}-1}\mathbh{1}_{\{
B_i\geq\delta n \}}\geq\frac{\alpha_0 K_0}{2}\Biggr)
\\[-3pt]
&&\qquad\quad{}+\PP\Biggl(\rho^{(1-\delta)n}\sum_{i=0}^{\nu_{x_n}-1}\mathbh
{1}_{\{B_i\geq\delta n \}}< \frac{\alpha_0
K_0}{2}\Biggr).
\nonumber
\end{eqnarray}
By the definition of $x_n$ and (\ref{eqslln}), $\nu_{x_n}-1$ is
asymptotically equivalent to~$(1-\allowbreak\rho)\lfloor\rho^{-n}t\rfloor$,
hence, if $\alpha_0$ is taken to be $(1-\rho)t$, by (\ref
{eq1}), one gets that the last
expression converges to $0$ as $n$ gets large. In the second term, the
mouse stays below level~$\delta n$, so a visit of the cat to $\delta n$
on a cycle is necessarily at least one meeting of the cat and the mouse
on this cycle. Further, it is clear that $\nu_{x_n}-1$ is not larger
than $x_n=\lfloor\rho^{-n}t\rfloor$. Finally, recall that the mouse
moves only when met by the cat and the sequence of successive
sites visited by the mouse is a also a simple reflected random walk.
Hence, if $\alpha_1=\alpha_0K_0/2$,
\begin{eqnarray*}
&&\PP\Biggl(\sup_{1\leq k\leq\lfloor\rho^{-n}t\rfloor}M_k\leq
\delta n,
\rho^{(1-\delta)n}\sum_{i=0}^{\nu_{x_n}-1}\mathbh{1}_{\{
B_i\geq\delta n \}}\geq
\alpha_1\Biggr)
\\[-3pt]
&&\qquad\leq\PP\Biggl(\sup_{1\leq k\leq\lfloor\rho^{-n}t\rfloor
}M_k\leq
\delta n,
\rho^{(1-\delta)n}\sum_{i=0}^{\lfloor\rho^{-n}t\rfloor}\mathbh
{1}_{\{C_i=M_i \}}\geq
\alpha_1\Biggr)
\\[-3pt]
&&\qquad\leq\PP\Bigl(\sup_{1\leq k\leq\lfloor\alpha_1\rho^{-(1-\delta
)n}\rfloor}C_k\leq\delta
n\Bigr)
=\PP\bigl(T_{\lfloor\delta n\rfloor+1}\geq\bigl\lfloor\alpha_1\rho
^{-(1-\delta)n}\bigr\rfloor\bigr)
\end{eqnarray*}
with the notation of Proposition \ref{ConvExp}, but this proposition
shows that the
random variable $\rho^{\lfloor\delta n\rfloor}T_{\lfloor\delta
n\rfloor+1}$ converges in
distribution as $n$ gets large. Consequently, since $\delta<1/2$, the
expression
\[
\PP\bigl(T_{\lfloor\delta n\rfloor+1}\geq\bigl\lfloor\alpha_1\rho
^{-(1-\delta)n}\bigr\rfloor\bigr)=
\PP\bigl(\rho^{\lfloor\delta n\rfloor}T_{\lfloor\delta n\rfloor
+1}\geq\alpha_1\rho^{-(1-2\delta)n}\bigr)
\]
converges to $0$. The relation
\[
\lim_{n\to+\infty} \PP\biggl(\sup_{0\leq u\leq t} \frac
{M_{\lfloor u\rho^{-n}\rfloor}}{n}\geq\frac{1}{2}\biggr)=1
\]
has been proved.\vadjust{\goodbreak}

The proof of the same result on the interval $[s,t]$ uses a coupling argument.
Define the cat and mouse Markov chain $(\widetilde{C}_k,\widetilde
{M}_k)$ as follows:
\[
(\widetilde{C}_k, k\geq0)=\bigl(C_{\lfloor s\rho^{-n}\rfloor+k}, k\geq0\bigr)
\]
and
the respective jumps of the sequences $(M_{\lfloor s\rho^{-n}\rfloor
+k})$ and $(\widetilde{M}_k)$ are independent except when
$M_{\lfloor s\rho^{-n}\rfloor+k}=\widetilde{M}_k$, in which\vspace*{1pt} case
they are the same. In
this way, one checks that $(\widetilde{C}_k,\widetilde{M}_k)$ is a
cat and mouse Markov
chain with the initial condition
\[
(\widetilde{C}_0,\widetilde{M}_0)=\bigl(C_{\lfloor
s\rho^{-n}\rfloor},0\bigr).
\]
By induction on $k$, one gets that
$M_{\lfloor s\rho^{-n}\rfloor+k}\geq\widetilde{M}_k$ for all
$k\geq0$.
Because of the ergodicity of $(C_k)$, the variable $C_{\lfloor s\rho
^{-n}\rfloor}$
converges in distribution as $n$ get large. Thus, $\widetilde{C}_0$ is
on a finite distance from 0 with probability one, and in the same way
as before, one gets that
\[
\lim_{n\to+\infty} \PP\biggl(\sup_{0\leq u\leq t-s} \frac
{\widetilde{M}_{\lfloor
u\rho^{-n}\rfloor}}{n}\geq\frac{1}{2}\biggr)=1,
\]
therefore,
\[
\liminf_{n\to+\infty} \PP\biggl(\sup_{s\leq u\leq t} \frac
{{M}_{\lfloor u\rho^{-n}\rfloor}}{n}\geq\frac{1}{2}\biggr)\geq
\liminf_{n\to+\infty} \PP\biggl(\sup_{0\leq u\leq t-s} \frac
{\widetilde{M}_{\lfloor u\rho^{-n}\rfloor}}{n}\geq\frac
{1}{2}\biggr)=1.
\]
This completes the proof of relation (\ref{Irancy}).
\end{pf*}

\section{Continuous time Markov chains}\label{secMMI}
Let $Q=(q(x,y), x, y\in{\cal S})$ be the $Q$-matrix of a continuous
time Markov chain on
${\cal S}$ such that, for any $x\in{\cal S}$,
\[
q_x\stackrel{\mathrm{def.}}{=}\sum_{y\dvtx y\not=x}q(x,y)
\]
is finite and that the Markov chain is positive recurrent and $\pi$ is
its invariant
probability distribution. The transition matrix of the
underlying discrete time Markov chain is denoted as
$p(x,y)=q(x,y)/q_x$; for
$x\not=y$, note that $p(\cdot,\cdot)$ vanishes on the diagonal.
See Norris \cite{Norris} for an introduction on Markov chains and
Rogers and
Williams \cite{Rogers} for a more advanced presentation.

The analogue of the Markov chain $(C_n,M_n)$ in this setting is the Markov
chain $(C(t),M(t))$ on ${\cal S}^2$ whose infinitesimal generator
$\Omega$ is defined by, for~$x$, $y\in{\cal S}$,
%
%
\begin{eqnarray}
\Omega(f)(x,y)&=& \sum_{z\in{\cal S}} q(x,z)[f(z,y)-f(x,y)]\mathbh
{1}_{\{x\not=y\}}
\nonumber\\[-8pt]\\[-8pt]
&&{}+ \sum_{z, z'\in{\cal S}} q_xp(x,z)p(x,z')[f(z,z')-f(x,x)]\mathbh
{1}_{\{x=y\}}
\nonumber
\end{eqnarray}
for any function $f$ on ${\cal S}^2$ vanishing outside a finite set.
The first coordinate
is indeed a Markov chain with $Q$-matrix $Q$ and when the cat and the
mouse are at the
same site $x$, after an exponential random time with parameter~$q_x$,
they jump
independently according to the transition matrix $P$. Note that if one
looks at the
sequence of sites visited by $(C(t),M(t))$, then it has the same
distribution as the cat
and mouse Markov chain associated to the matrix $P$. For this reason,
the results
obtained in Section \ref{secCM} can be proved easily in this setting.
In particular,
$(C(t),M(t))$ is null recurrent when $(C(t))$ is reversible.
\begin{prop}
\label{propmst}
If, for $t\geq0$,
\[
U(t)=\int_0^{t}\mathbh{1}_{\{M(s)=C(s)\}} \,ds
\]
and $S(t)=\inf\{s>0\dvtx U(s)\geq t\}$, then the process $(M(S(t)))$ has
the same distribution
as $(C(t))$, that is, it is a Markov process with $Q$-matrix $Q$.
\end{prop}

This proposition simply states that, up to a time change, the mouse
moves like the cat. In discrete time this is fairly obvious; the proof
is in this case a~little more technical.
\begin{pf*}{Proof of Proposition \ref{propmst}}
If $f$ is a function on ${\cal S}$, then by characterization of Markov
processes, one has that the process
\[
(H(t))\stackrel{\mathrm{def.}}{=}\biggl(f(M(t))-f(M(0))-\int_0^t
\Omega(\bar{f})(C(s),M(s)) \,ds\biggr)
\]
is a local martingale with respect\vspace*{1pt} to the natural filtration $({\cal
F}_t)$ of
$(C(t),M(t))$,
where $\bar{f}\dvtx{\cal S}^2\to\R$ such that $\bar{f}(x,y)=f(y)$ for
$x$, $y\in{\cal
S}$. The fact that, for $t\geq0$, $S(t)$ is a stopping time and that
$s\to S(s)$ is nondecreasing, and
Doob's optional stopping theorem imply that $(H(S(t)))$ is a local martingale
with respect to the filtration $({\cal F}_{S(t)})$. Since
\begin{eqnarray*}
&&\int_0^{S(t)} \Omega(\bar{f})(C(s),M(s)) \,ds\\
&&\qquad=\sum_{y\in{\cal S}} \int_0^{S(t)}q(M(s),y)\mathbh{1}_{\{
C(s)=M(s)\}}
\bigl(f(y)-f(M(s))\bigr) \,ds\\
&&\qquad=\int_0^{S(t)}\mathbh{1}_{\{C(s)=M(s)\}} Q(f)(M(s))\,
ds\\
&&\qquad= \int_0^{t} Q(f)(M(S(s))) \,ds,
\end{eqnarray*}
the infinitesimal generator $Q$ is defined for $x\in{\cal S}$ in a
standard way as
\[
Q(f)(x)=\sum_{y\in{\cal S}}q(x,y)[f(y)-f(x)].
\]
One therefore gets that
\[
\biggl(f(M(S(t)))-f(M(0))-\int_0^t Q(f)(M(S(s))) \,ds\biggr)
\]
is a local martingale for any function $f$ on ${\cal S}$. This implies
that $(M(S(t)))$ is
a Markov process with $Q$-matrix $Q$, that is, that $(M(S(t)))$ has the
same distribution as
$(C(t))$. See Rogers and Williams \cite{Rogers06}.
\end{pf*}

\subsection*{The example of the $M/M/\infty$ process}
The example of the $M/M/\infty$ queue is investigated in the rest of
this section.
The associated Markov process can be seen as an example of a discrete
Ornstein--Uhlenbeck process.
As it will be shown, there is a significant qualitative difference with the
example of Section \ref{secMM1} which is a discrete time version of
the $M/M/1$
queue. The $Q$-matrix is given by
%
%
\begin{equation}\label{OU}
\cases{q(x,x+1)= \rho, \cr
q(x,x-1) = x.}
\end{equation}
The corresponding Markov chain is positive recurrent and reversible and
its invariant
probability distribution is Poisson with parameter $\rho$.
\begin{prop}\label{Hit+}
If $C(0)=x\leq n-1$ and
\[
T_n=\inf\{s>0\dvtx C(s)=n\},
\]
then, as $n$ tends to infinity, the variable $T_n/E_x(T_n)$ converges
in distribution to an
exponentially distributed random variable with parameter $1$ and
\[
\lim_{n\to+\infty} \E_x(T_n){\rho^n}/{(n-1)!}=e^{-\rho}.
\]
If $C(0)=n$, then $T_0/\log n$ converges in distribution to $1$.
\end{prop}

See Chapter 6 of Robert \cite{Robert08}. It should be
remarked that the
duration of time it takes to reach $n$ starting from $0$ is essentially
the time it takes
to go to $n$ starting from $n-1$.

\subsection*{Multiplicative jumps}
The above proposition gives the order of magnitude for the duration of
time for the cat to hit the mouse. As before, the cat returns
``quickly'' to the neighborhood of $0$, but, contrary to the reflected
random walk, it turns out that the cat will take the mouse down for
some time before leaving the mouse. The next proposition shows that if
the mouse is at $n$, its next location after the visit of the cat is
of the order of $nF$ for a~certain random variable $F$.
\begin{prop}\label{MMIScal}
If $C(0)=M(0)=n$ and
\[
T_0=\inf\{s>0\dvtx C(s)=0\},
\]
then, as $n$ goes to infinity, the random variable $M(T_0)/n$ converges
in distribution to
a random variable $F$ on $[0,1]$ such that $\PP(F\leq x)=x^\rho$.
\end{prop}
\begin{pf}
Let $\tau=\inf\{s>0\dvtx M(s)=M(s-)+1\}$ be the instant of the first
upward jump of
$(M(s))$. Since $(M(S(s)))$ has the same distribution as~$(C(s))$, one
gets that $U(\tau)$, with $U(t)$ defined as in Proposition \ref
{propmst}, has the same distribution as the time till a first upward
jump of $(C(s))$, which is an exponential random variable with
parameter $\rho$ by definition (\ref{OU}). Now, think of $(M(S(s)))$
as the process describing the number of customers in an $M/M/\infty$
queue, which contains $n$ customers at time $0$. Let $(E_i)$ be i.i.d.
exponential random variables with parameter~$1$. For $1\leq i\leq n$,
$E_i$ is the service time of the $i$th initial customer. At time $\tau
$, the process of
the mouse will have run only for $U(\tau)$, so the $i$th customer is
still there if $E_i> U(\tau)$. Note that there is no arrival up
to time $\tau$, and, hence,
\[
M(\tau)\stackrel{\mathrm{dist.}}{=}1+\sum_{i=1}^n\mathbh{1}_{
\{E_i> U(\tau)\}}.
\]
Consequently, by conditioning on the value of $U(\tau)$, by the law of
large numbers
one obtains that the sequence $(M(\tau)/n)$ converges in distribution
to the random variable
$F\stackrel{\mathrm{def.}}{=}\exp(-U(\tau))$, which implies
directly that $\PP(F\leq x)=x^\rho$.

It remains to show that $(M(T_0)/n)$ converges in distribution to the
same limit as $(M(\tau)/n)$. The fact that the mouse moves only when
it meets the cat gives the following:
\begin{itemize}[--]
\item[--] On the event $\tau\geq T_0$, necessarily $M(\tau-)=C(\tau-)=0$
because if the mouse did not move upward before time $T_0$, then it has
reached $0$ together with the cat. In this case, at time $\tau$, the
mouse makes its first jump upward from $0$ to 1. Thus, the quantity
\[
\PP(\tau\geq T_0)\leq\PP\bigl(M(\tau)=1\bigr)=\PP\bigl(U(\tau)>\max\{
E_1,\ldots,
E_n\}\bigr)
\]
converges to $0$ as $n\to\infty$.
\item[--] Just before time $\tau$, the mouse and the cat are at the same
location and
\[
\PP\bigl(C(\tau)=M(\tau-)-1\bigr)=\E\biggl[\frac{M(\tau-)}{\rho+M(\tau
-)}\biggr]
\]
converges to $1$ as $n$ gets large.
\end{itemize}
The above statements imply that with probability converging to one, the
cat will find itself below the mouse for the\vadjust{\goodbreak} first time strictly above
level zero and before time $T_0$. We now show that after this event the
cat will hit zero before returning back to the mouse. If $\eps>0$, then
\begin{eqnarray*}
\E\bigl(\PP_{C(\tau)}\bigl(T_0\geq T_{M(\tau)}\bigr)\bigr)&\leq&
\E\bigl(\PP_{M(\tau)-1}\bigl(T_0\geq T_{M(\tau)}\bigr)\bigr)\\
&\le&\PP\biggl(\frac{M(\tau)}{n}\leq\eps\biggr)+\sup_{k\geq
\lfloor\eps n\rfloor} \PP_k(T_0\geq T_{k+1}),
\end{eqnarray*}
hence, by Proposition \ref{Hit+}, for $\eps$ (resp., $n$) sufficiently small
(resp., large), the above quantity is arbitrarily small. This result
implies that
the probability of the event $\{M(\tau)=M(T_0)\}$
converges to $1$. The proposition is proved.
\end{pf}

\subsection*{An underlying random walk}
If $C(0)=0$ and $M(0)=n$, the next time the cat returns to $0$,
Proposition \ref{MMIScal}
shows that the mouse will be at a location of the order of $n F_1$, where
$F_1=\exp(-E_1/\rho)$ and $E_1$
is an exponential random variable with parameter $1$. After the $p$th round,
the location of the mouse is of the order of
%
%
\begin{equation}\label{Imul}
n\prod_{k=1}^pF_k=n\exp\Biggl(-\frac{1}{\rho}\sum_{k=1}^pE_k\Biggr),
\end{equation}
where $(E_k)$ are i.i.d. with the same distribution as $E_1$. A precise
statement of this
nonrigorous statement can be formulated easily. From (\ref
{Imul}), one gets that
after the order of $\rho\log n$ rounds, the location of the mouse is
within a~finite interval.

The corresponding result for the reflected random walk exhibits an
additive behavior.
Theorem \ref{TheoMM1} gives that the location of the mouse is of the
order of
%
%
\begin{equation}\label{Iadd}
n+\sum_{i=1}^p A_i
\end{equation}
after $p$ rounds, where $(A_k)$ are i.i.d. copies of $M'_\infty$,
distribution of which is given by the generating
function of relation (\ref{GenD}). In this case the number of rounds
after which the
location of the mouse is located within a finite interval is of the
order of~$n$.

As Theorem \ref{TheoMM1} shows, for the reflected random walk, $t\to
\rho^{-n}t$ is a
convenient time scaling to describe the location of the mouse until it
reaches a~finite
interval. This is not the case for the $M/M/\infty$ queue, since the
duration of the
first round of the cat, of the order of ${(n-1)!}/{{\rho^n}}$ by
Proposition~\ref{Hit+},
dominates by far the duration of the subsequent rounds, that is, when
the location of the
mouse is at $xn$ with $x<1$.


%

%
\printaddresses

\end{document}